\documentclass[12pt]{smfart}

\usepackage{smfthm}
\usepackage{smfenum}
\usepackage{amssymb}
\usepackage[french]{babel}

\newcommand{\on}{\operatorname}

\renewcommand{\tilde}{\widetilde}
\newcommand{\Qp}{\mathbf{Q}_p}
\newcommand{\Zp}{\mathbf{Z}_p}
\newcommand{\Cp}{\mathbf{C}_p}
\newcommand{\ZZ}{\mathbf{Z}}
\newcommand{\RR}{\mathbf{R}}
\newcommand{\QQ}{\mathbf{Q}}
\newcommand{\OO}{\mathcal{O}}
\newcommand{\amax}{\mathbf{A}_{\mathrm{max}}}
\newcommand{\atplus}{\widetilde{\mathbf{A}}^+}
\newcommand{\atdag}[1]{\widetilde{\mathbf{A}}^{\dagger #1}}
\newcommand{\adag}[1]{\mathbf{A}^{\dagger #1}}
\newcommand{\btrig}[2]{\widetilde{\mathbf{B}}^{\dagger #1}_{\mathrm{rig} #2}}
\newcommand{\bnrig}[2]{\mathbf{B}^{\dagger #1}_{\mathrm{rig} #2}}
\newcommand{\btrigplus}[1]{\widetilde{\mathbf{B}}^{+}_{\mathrm{rig} #1}}
\newcommand{\bcris}{\mathbf{B}_{\mathrm{cris}}}
\newcommand{\bst}{\mathbf{B}_{\mathrm{st}}}
\newcommand{\bdR}{\mathbf{B}_{\mathrm{dR}}}
\newcommand{\btdag}[1]{\widetilde{\mathbf{B}}^{\dagger #1}}
\newcommand{\bdag}[1]{\mathbf{B}^{\dagger #1}}
\newcommand{\bb}{\mathbf{B}}
\newcommand{\et}{\tilde{\mathbf{E}}}
\newcommand{\dnrig}[1]{\mathbf{D}^{\dagger #1}_{\mathrm{rig}}}
\newcommand{\dst}{\mathbf{D}_{\mathrm{st}}}
\newcommand{\dstp}[1]{\mathbf{D}_{\mathrm{st},#1}}
\newcommand{\ddR}{\mathbf{D}_{\mathrm{dR}}}
\renewcommand{\ddag}[1]{\mathbf{D}^{\dagger #1}}
\newcommand{\ndr}{\mathbf{N}_{\mathrm{dR}}}
\newcommand{\eps}{\varepsilon}
\newcommand{\ra}{\rightarrow}
\newcommand{\lra}{\longrightarrow}
\newcommand{\pibar}{\overline{\pi}}
\renewcommand{\phi}{\varphi}

\newcommand{\logpi}{\log[\overline{\pi}]}
\newcommand{\mbf}{\mathbf{M}}
\newcommand{\pgm}{\mathcal{M}}
\newcommand{\nbf}{\mathbf{N}}
\newcommand{\dbf}{\mathbf{D}}
\newcommand{\fil}{\mathrm{Fil}}
\renewcommand{\geq}{\geqslant}
\renewcommand{\leq}{\leqslant} 

\title{{\'E}quations diff{\'e}rentielles $p$-adiques et $(\phi,N)$-modules
filtr{\'e}s}

\alttitle{$p$-adic differential equations and filtered
  $(\phi,N)$-modules} 

\author{Laurent Berger}

\address{52 Rue de Nanterre \\ 92600 Asni{\`e}res \\ France}

\email{laurent@math.harvard.edu}

\urladdr{www.math.harvard.edu/\~{}laurent}

\subjclass{11F80, 12H25, 13K05, 14F30.}

\keywords{repr{\'e}sentations $p$-adiques, de de Rham,
  $(\phi,N)$-modules filtr{\'e}s, {\'e}quations diff{\'e}rentielles
  $p$-adiques} 

\date{Juin 2004}

\NumberTheoremsIn{subsection}

\begin{document}

\begin{abstract}
L'objet de cet article est de montrer que les deux cat{\'e}gories
suivantes sont {\'e}quivalentes (1) la cat{\'e}gorie des
$(\phi,N,G_K)$-modules filtr{\'e}s (2) la cat{\'e}gorie des
$(\phi,\Gamma_K)$-modules sur l'anneau de Robba tels que l'alg{\`e}bre
de Lie de $\Gamma_K$ agit localement trivialement.

De plus, on montre que sous cette {\'e}quivalence, les
$(\phi,N,G_K)$-modules filtr{\'e}s admissibles correspondent aux
$(\phi,\Gamma_K)$-modules {\'e}tales, ce qui nous permet de donner une
nouvelle d{\'e}monstration du th{\'e}or{\`e}me de Colmez-Fontaine. 
\end{abstract}

\begin{altabstract}
The goal of this article is to show that the following two categories
are equivalent (1) the category of filtered $(\phi,N,G_K)$-modules (2)
the category of $(\phi,\Gamma_K)$-modules over the Robba ring such
that the Lie algebra of $\Gamma_K$ acts locally trivially.

Furthermore, we show that under this equivalence, the admissible 
filtered $(\phi,N,G_K)$-modules correspond to the {\'e}tale
$(\phi,\Gamma_K)$-modules, which gives a new proof of
Colmez-Fontaine's theorem.
\end{altabstract}

\maketitle
\setcounter{tocdepth}{2}
\tableofcontents
\setlength{\baselineskip}{18pt}

\section*{Introduction}
Dans tout cet article, $K$ est un corps local 
qui contient $\Qp$ et qui est
muni d'une valuation discr{\`e}te {\'e}tendant la
valuation $p$-adique et pour laquelle $K$ est complet 
et de corps r{\'e}siduel parfait $k_K$. Pour $n \leq \infty$, on pose
$K_n=K(\mu_{p^n})$ et $\Gamma_K=\on{Gal}(K_\infty/K)$. 

\subsection*{$(\phi,N,G_K)$-modules filtr{\'e}s et $(\phi,\Gamma_K)$-modules}
L'objet de cet article est de cons\-truire une {\'e}quivalence de
cat{\'e}gories entre deux cat{\'e}gories utilis{\'e}es dans la th{\'e}orie des
repr{\'e}sentations $p$-adiques et des {\'e}quations diff{\'e}rentielles
$p$-adiques (on se reportera au chapitre \ref{secI} pour des rappels
sur ces cat{\'e}gories). Il s'agit de:
\begin{enumerate}
\item la cat{\'e}gorie des $(\phi,N,G_K)$-modules
filtr{\'e}s, dont la sous-cat{\'e}gorie des objets admissibles
param{\'e}trise les repr{\'e}\-sentations
$p$-adiques potentiellement semi-stables du groupe de Galois absolu
$G_K$ du corps $K$;
\item  la cat{\'e}gorie des
$(\phi,\Gamma_K)$-modules sur l'anneau de Robba tels que l'alg{\`e}bre
de Lie de $\Gamma_K$ agit localement trivialement,
qui g{\'e}n{\'e}ralise
la notion d'{\'e}quation diff{\'e}rentielle $p$-adique munie d'une structure
de Frobenius. 
\end{enumerate}

On construit un foncteur $D \mapsto
\pgm(D)$ qui {\`a} un $(\phi,N,G_K)$-module filtr{\'e} $D$ associe un
$(\phi,\Gamma_K)$-module $\pgm(D)$ sur l'anneau de Robba
$\bnrig{}{,K}$ et le r{\'e}sultat principal de cet article est le suivant:

\begin{enonce*}{Th{\'e}or{\`e}me A}
Le $\otimes$-foncteur exact
$D \mapsto \pgm(D)$, de la cat{\'e}gorie des
$(\phi,N,G_K)$-modules filtr{\'e}s, dans la cat{\'e}gorie des
$(\phi,\Gamma_K)$-modules sur $\bnrig{}{,K}$ dont la connexion
associ{\'e}e est localement triviale, est une {\'e}quivalence de
cat{\'e}gories.
\end{enonce*}

{\'E}tant donn{\'e} un $(\phi,N,G_K)$-module filtr{\'e} $D$, le
$(\phi,\Gamma_K)$-module $\pgm(D)$ admet 
(par le th{\'e}or{\`e}me \cite[th{\'e}or{\`e}me 6.10]{KK00} de Kedlaya)
une filtration canonique par les 
{\og pentes de Frobenius \fg}. Nous calculons
ces pentes en terme des invariants $t_H$ et $t_N$ de $D$ et de ses
sous-objets. Comme corollaire de ces calculs, nous obtenons le
r{\'e}sultat suivant:

\begin{enonce*}{Th{\'e}or{\`e}me B}
Le $(\phi,N,G_K)$-module filtr{\'e} $D$ est admissible si et seulement
si $\pgm(D)$ est {\'e}tale. 
\end{enonce*}

\subsection*{Application aux repr{\'e}sentations $p$-adiques} 
Comme application du th{\'e}or{\`e}me B, on obtient une nouvelle
d{\'e}monstration du th{\'e}or{\`e}me \cite[th{\'e}or{\`e}me A]{CF00} 
de Colmez-Fontaine: 

\begin{enonce*}{Th{\'e}or{\`e}me}
Si $D$ est un
$(\phi,N,G_K)$-module filtr{\'e} admissible, alors il existe une
repr{\'e}sen\-tation $p$-adique $V$ potentiellement semi-stable telle que
$\dbf_{\rm pst}(V)=D$.
\end{enonce*}

Enfin, les constructions pr{\'e}c{\'e}dentes nous permettent de pr{\'e}ciser
les liens entre les divers invariants que l'on peut associer {\`a} une
repr{\'e}sentation $p$-adique potentiellement semi-stable $V$ (qui
devient semi-stable sur une extension galoisienne finie $L$ de $K$).  
Ces invariants sont: 
\begin{enumerate}
\item le $(\phi,N,G_{L/K})$-module filtr{\'e} $\dstp{L}(V)=(\bst
  \otimes_{\Qp} V)^{G_L}$;
\item l'{\'e}quation diff{\'e}rentielle $p$-adique $\ndr(V)$;
\item le $(\phi,\Gamma_K)$-module
{\'e}tale sur l'anneau de Robba $\dnrig{}(V)$.
\end{enumerate}

{\`A} l'{\'e}quation diff{\'e}rentielle $p$-adique $\ndr(V)$, on associe (cf
d{\'e}finition \ref{defisol}) son espace $\on{Sol}_L(\ndr(V))$ de
$G_L$-solutions, qui est un $(\phi,N,G_{L/K})$-module,
et on montre dans le th{\'e}or{\`e}me \ref{imageessent}
comment la donn{\'e}e de $\dnrig{}(V)$ d{\'e}termine une filtration sur
$\on{Sol}_L(\ndr(V))$, ce qui en fait un $(\phi,N,G_{L/K})$-module
filtr{\'e}. On a alors (les notations restantes 
sont d{\'e}finies dans le corps du texte): 

\begin{enonce*}{Th{\'e}or{\`e}me C}
Si $V$ est une repr{\'e}sentation $p$-adique de $G_K$ semi-stable quand
on la restreint {\`a} $G_L$, alors:
\begin{enumerate}
\item $\dstp{L}(V) = \on{Sol}_L(\ndr(V))$; 
\item $\dnrig{}(V)=\pgm(\dstp{L}(V))$;
\item $\ndr(V)=(\bnrig{}{,L}[\ell_X] \otimes_{L_0} 
\dstp{L}(V))^{\on{Gal}(L_\infty/K_\infty),N=0}$.
\end{enumerate}
Ces identifications sont compatibles {\`a} toutes les structures en pr{\'e}sence.
\end{enonce*}

Cela permet notamment de {\og retrouver la filtration \fg} sur
$\dstp{L}(V)$ quand on le construit {\`a} partir de $\dnrig{}(V)$
(cf \cite[proposition 5.6]{PC01} pour une construction similaire).

Enfin, nous montrons comment 
reconstruire, pour une repr{\'e}sentation semi-stable $V$,
le $(\phi,\Gamma_K)$-module $\ddag{}(V)$ sur $\bdag{}_K$
{\`a} partir de $\dst(V)$ (comme dans \cite[\S 4.3.1]{PC04}). 

\begin{enonce*}{Th{\'e}or{\`e}me D}
Si $V$ est une repr{\'e}sentation semi-stable de $G_K$ et 
si 
\begin{itemize}
\item $\{e_i\}_{i=1 \cdots d}$ est une base de $\dst(V)$ adapt{\'e}e
{\`a} la d{\'e}composition par les pentes de Frobenius;
\item $N(e_i)=\sum_{j=1}^d n_{j,i} e_j$;
\item $\{f_j\}_{j=1 \cdots d}$ est une base de $\ddR(V)$ adapt{\'e}e {\`a} la filtration
et $\phi^{-n}(e_i) = \sum_{j=1}^d p_{j,i}^{(n)} f_j$, 
\end{itemize}
alors $x = \sum_{i=1}^d x_i(X) \otimes e_i 
\in \bnrig{,r}{,K}[\ell_X,1/t] \otimes_{K_0} \dst(V)$ appartient {\`a}
$\ddag{,r}(V)$ si et seulement si:
\begin{enumerate}
\item $N(x_j)+\sum_{i=1}^d n_{j,i} x_i = 0$ pour $j=1 \cdots d$;
\item $\sum_{i=1}^d \iota_n(x_i)p_{j,i}^{(n)} \in t^{-t_H(f_j)}
  K_n[\![t]\!]$ pour $j=1 \cdots d$ et $n \geq n(r)$;
\item $\on{ord}(x_i) \leq - \on{pente}(e_i)$ pour $i=1 \cdots d$.
\end{enumerate}
\end{enonce*}

Ceci permet de construire {\og explicitement \fg} les {\'e}l{\'e}ments de
$\ddag{,r}(V)$ ce qui a des applications directes {\`a} la th{\'e}orie des 
repr{\'e}sentations de $\on{GL}(2,\Qp)$ de Breuil. 

\vspace{\baselineskip}
\noindent\textbf{Remerciements:} Je remercie Pierre Colmez, Jean-Marc
Fontaine et Kiran Kedlaya pour des discussions {\'e}clairantes sur certains
points de cet article. 

\renewcommand{\thesection}{\Roman{section}}
\section{Rappels et compl{\'e}ments}\label{secI}
Dans tout cet article, $K$ est un corps local 
qui contient $\Qp$ et qui est
muni d'une valuation discr{\`e}te {\'e}tendant la
valuation $p$-adique et pour laquelle $K$ est complet 
et de corps r{\'e}siduel parfait $k_K$. 

On {\'e}crit $\mu_{p^n} \subset \overline{K}$ 
(la cl{\^o}ture alg{\'e}brique de $K$) pour d{\'e}signer 
l'ensemble des racines $p^n$-i{\`e}mes de l'unit{\'e}, et 
pour $n \geq 1$, on d{\'e}finit
$K_n=K(\mu_{p^n})$ ainsi que $K_{\infty} = \cup_{n=1}^{+\infty}
K_n$. Si $n=0$, on pose $K_0 = W(k_K)[1/p]$ 
ce qui fait que $K/K_0$ est totalement ramifi{\'e}e. 
Les notations ne sont donc pas
vraiment compatibles, mais elles sont usuelles. 
Finalement, $K_0'$ d{\'e}signe 
l'extension maximale non-ramifi{\'e}e de $K_0$
contenue dans $K_\infty$. 

Soient $G_K=\on{Gal}(\overline{K}/K)$ et 
$H_K=\on{Gal}(\overline{K}/K_\infty)$ 
le noyau du caract{\`e}re cyclotomique 
$\chi: G_K \ra \Zp^*$ et $\Gamma_K=G_K/H_K$
le groupe de Galois de $K_{\infty}/K$,
qui s'identifie via le caract{\`e}re 
cyclotomique {\`a} un sous-groupe ouvert de 
$\Zp^*$. On note $\sigma$ le Frobenius absolu 
(qui rel{\`e}ve $x \mapsto
x^p$ sur $k_K^{\on{sep}}$).

\subsection{Les $(\phi,N,G_{L/K})$-modules filtr{\'e}s}\label{pnfildef}
Le corps $L$ d{\'e}signe une extension galoisienne finie de $K$, et
$G_{L/K}$ le groupe de Galois de $L/K$.
La cat{\'e}gorie des $(\phi,N,G_{L/K})$-modules filtr{\'e}s et sa 
sous-cat{\'e}gorie pleine de $(\phi,N,G_{L/K})$-modules filtr{\'e}s admissibles
sont {\'e}tudi{\'e}es en d{\'e}tail dans \cite[\S 4]{Bu88sst}. Nous nous
contentons donc de rappeler les d{\'e}finitions et quelques r{\'e}sultats
dont nous aurons besoin par la suite. 

Un $(\phi,N,G_{L/K})$-module est un $L_0$-espace vectoriel $D$ de
dimension finie et muni d'une application $\sigma$-semi-lin{\'e}aire
bijective $\phi:D \ra D$, d'une application lin{\'e}aire $N:D \ra D$ 
qui v{\'e}rifie $N \phi = p \phi N$ et d'une action semi-lin{\'e}aire de
$G_{L/K}$  qui commute {\`a} $\phi$ et $N$. 
Si $L=K$, on parle tout simplement de $(\phi,N)$-module 
(relatif {\`a} $K$). 

Un $(\phi,N,G_{L/K})$-module filtr{\'e} est la donn{\'e}e 
d'un $(\phi,N,G_{L/K})$-module
$D$ et d'une filtration d{\'e}croissante exhaustive et separ{\'e}e
$\fil^i D_L$ sur $D_L = L \otimes_{L_0} D$, par des sous $L$-espaces
vectoriels stables par $G_{L/K}$. Il est {\'e}quivalent de se donner une
filtration sur $D_K = D_L^{G_{L/K}}$. 

La cat{\'e}gorie des $(\phi,N,G_{L/K})$-modules filtr{\'e}s est une
cat{\'e}gorie additive $\Qp$-lin{\'e}aire qui admet des noyaux et des
conoyaux (mais qui n'est pas ab{\'e}lienne), ainsi que des produits
tensoriels et des $\on{Hom}$ internes. 

Par abus de langage, on dira que $D$ est un $(\phi,N,G_K)$-module
filtr{\'e} s'il existe une extension finie galoisienne $L/K$ telle que
$D$ est un $(\phi,N,G_{L/K})$-module filtr{\'e}. 

Si $D$ est un $(\phi,N,G_{L/K})$-module filtr{\'e} de dimension $1$, on
d{\'e}finit $t_H(D)$ comme le plus grand entier $i$ tel que $\fil^i D_L
\neq 0$ et si $\phi(d)=\lambda d$ avec $d \in D$, alors
$v_p(\lambda)$ ne d{\'e}pend pas du choix de $d \neq 0$ et on 
d{\'e}finit $t_N(D)=v_p(\lambda)$. 
Si $D$ est un $(\phi,N,G_{L/K})$-module
filtr{\'e} de dimension $\geq 1$, on d{\'e}finit $t_H(D)=t_H(\det D)$
et $t_N(D)=t_N(\det D)$. Si $e \in D_L$, on pose $t_H(e)=t_H(L \cdot
e)$. 

On dit qu'un $(\phi,N,G_{L/K})$-module filtr{\'e} $D$ est admissible si
$t_H(D)=t_N(D)$ et si pour tout sous-objet $D'$ de $D$, on a
$t_N(D')-t_H(D') \geq 0$. La cat{\'e}gorie des
$(\phi,N,G_{L/K})$-modules filtr{\'e}s admissibles est une
sous-cat{\'e}gorie pleine de la cat{\'e}gorie des
$(\phi,N,G_{L/K})$-modules filtr{\'e}s et elle est de plus ab{\'e}lienne. 

Rappelons que la principale source de $(\phi,N,G_{L/K})$-modules
filtr{\'e}s admissibles est la suivante: si $V$ est une repr{\'e}sentation
$p$-adique de $G_K$ dont la restriction {\`a} $G_L$ est semi-stable, alors
$\dstp{L}(V)$ est un $(\phi,N,G_{L/K})$-module
filtr{\'e} admissible.  

\subsection{L'anneau de Robba}\label{robbarap}
D{\'e}finissons ici quelques anneaux de s{\'e}ries formelles
(ces constructions sont faites en d{\'e}tail dans \cite{PC03}). 
Si $r$ est
un r{\'e}el positif et $F=K_0$ (pour all{\'e}ger un peu les notations), 
soit $\bdag{,r}_F$ l'anneau des s{\'e}ries formelles
$f(X)=\sum_{k \in \ZZ} a_k X^k$ o{\`u} $\{a_k \in F\}_{k \in \ZZ}$ est une suite
born{\'e}e telle que $f(X)$ converge sur la couronne $0 < v_p(X) \leq
1/r$. Cet anneau est muni
d'une action de $\Gamma_F$, qui est triviale sur les coefficients
et donn{\'e}e par $\gamma(X)=(1+X)^{\chi(\gamma)}-1$ et on peut
d{\'e}finir un Frobenius $\phi: \bdag{,r}_F \ra \bdag{,pr}_F$ qui est
$\sigma$-semi-lin{\'e}aire sur les coefficients
et tel que $\phi(X)=(1+X)^p-1$. Le {\og th{\'e}or{\`e}me de 
pr{\'e}paration de Weierstrass \fg} montre que $\bdag{}_F = \cup_{r \geq
0} \bdag{,r}_F$ est un corps. Ce corps n'est pas complet pour la
norme de Gauss et on appelle $\bb_F$ son
compl{\'e}t{\'e} qui est un corps local de dimension $2$ dont le corps
r{\'e}siduel s'identifie {\`a} $k_K(\!(\overline{X})\!)$.
 
L'extension $K_\infty / F_\infty$ est une extension finie de 
degr{\'e} de ramification $e_K \leq [K_\infty:F_\infty]$ 
et par la th{\'e}orie du corps de
normes de \cite{FW79,WI83} il lui correspond une extension s{\'e}parable
$k_K'(\!(\overline{X}_K)\!) / k_K(\!(\overline{X})\!)$ de 
degr{\'e} $[K_\infty:F_\infty]$ qui nous permet de d{\'e}finir des
extensions non-ramifi{\'e}es $\mathbf{B}_K / \mathbf{B}_F$ et 
$\bdag{}_K / \bdag{}_F$ de degr{\'e} $[K_\infty:F_\infty]$. 
On peut montrer que 
$\bdag{}_K = \cup_{r \geq 0} \bdag{,r}_K$ 
et qu'il existe $r_0(K)$ tel que si $r \geq r_0(K)$, alors   
$\bdag{,r}_K$ est un $\bdag{,r}_F$-module libre de rang
$[K_\infty:F_\infty]$ qui s'identifie
{\`a} un anneau de s{\'e}ries formelles
$f(X_K)=\sum_{k \in \ZZ} a_k X_K^k$ o{\`u} 
$\{a_k \in K_0'\}_{k \in \ZZ}$ est une suite
born{\'e}e telle que $f(X_K)$ converge sur la couronne 
$0 < v_p(X_K) \leq 1/e_K r$. L'{\'e}l{\'e}ment $\overline{X}_K$ 
v{\'e}rifie une {\'e}quation d'Eisenstein
sur $k_K'(\!(\overline{X})\!)$  
qu'on peut relever en une {\'e}quation sur
$\bdag{,r}_{K_0'}$; l'action de $\Gamma_K$
s'{\'e}tend naturellement {\`a} $\bdag{,r}_K$ de m{\^e}me que 
le Frobenius $\phi: \bdag{,r}_K \ra \bdag{,pr}_K$.

L'anneau $\bdag{,r}_K$ s'identifiant {\`a} un anneau de s{\'e}ries
formelles convergeant sur une couronne, 
il est naturellement muni d'une topologie de Fr{\'e}chet, 
la topologie de la {\og convergence compacte \fg}
sur les couronnes $C_K[r;s] = \{ z \in \OO_{\Cp},
\ 1/e_K s \leq v_p(z) \leq 1/e_K
r\}$ pour $s \geq r$, 
et son compl{\'e}t{\'e} $\bnrig{,r}{,K}$ 
pour cette topologie s'identifie {\`a} 
l'anneau de s{\'e}ries formelles
$f(X_K)=\sum_{k \in \ZZ} a_k X_K^k$ o{\`u} $\{a_k \in 
K_0'\}_{k \in \ZZ}$ est une suite
non n{\'e}cessairement born{\'e}e telle que 
$f(X_K)$ converge sur la couronne $0 < v_p(X_K) \leq 1/e_K r$.
Par exemple, si on pose $t=\log(1+X)$, alors 
$t \in \bnrig{,r}{,F} \subset 
\bnrig{,r}{,K}$ pour tout $r \geq 0$.
L'anneau $\bnrig{}{,K} = \cup_{r \geq 0} \bnrig{,r}{,K}$ est 
{\og l'anneau de Robba \fg}.  
Nous allons rappeler quelques-uns des r{\'e}sultats
de \cite[\S 4]{LB02} qui nous seront utiles dans la suite.

Il existe $r(K)$ que l'on peut supposer
$\geq r_0(K)$ tel que si $p^{n-1}(p-1) \geq 
r \geq r(K)$, alors on a une application injective $\iota_n :
\bdag{,r}_K \ra K_n[\![t]\!]$ (c'est l'application $\phi^{-n}$ de
\cite[\S III.2]{CC99}). 
Par exemple si $K=K_0$, alors
$\iota_n(X)=\eps^{(n)} \exp(t/p^n)-1$ o{\`u} $\eps^{(n)}$ est une racine
primitive $p^n$-i{\`e}me de $1$ et $\iota_n$ agit par $\sigma^{-n}$ sur les
coefficients. On d{\'e}finit $n(r)$ comme {\'e}tant le plus petit entier
$n$ tel que $p^{n-1}(p-1) \geq r$ ce qui fait que $\iota_n :
\bdag{,r}_K \ra K_n[\![t]\!]$ est d{\'e}finie d{\`e}s que $n \geq n(r)$. 

L'application $\iota_n$
se prolonge en une application injective 
$\iota_n : \bnrig{,r}{,K} \ra K_n[\![t]\!]$. 
L'action de $\Gamma_K$ sur $\bnrig{,r}{,K}$ s'{\'e}tend en une action de
l'alg{\`e}bre de Lie de $\Gamma_K$ donn{\'e}e par $\nabla(f) =
\log(\gamma)(f)/\log_p(\chi(\gamma))$ pour $\gamma \in \Gamma_K$ assez
proche de $1$. 
Si $f = f(X) \in \bnrig{,r}{,F}$ alors
$\nabla(f(X)) = t(1+X)df/dX$. Si $f \in \bnrig{,r}{,K}$ alors
on pose $\partial(f) = t^{-1} \nabla(f)$ ce qui fait que si 
$f = f(X) \in \bnrig{,r}{,F}$ alors
$\partial(f(X)) = (1+X)df/dX$ et que si $f \in \bnrig{}{,K}$ v{\'e}rifie
une {\'e}quation alg{\'e}brique $P(f)=0$ sur $\bnrig{}{,F}$ 
telle que $P'(f) \neq 0$, alors
on peut aussi calculer $\partial(f)$ par la formule
$\partial(f)= -(\partial P)(f)/P'(f)$. En particulier 
$\partial(f)=0$ si et seulement si $f \in K_0'$.

\begin{lemm}\label{partunit}
Si $w \geq 1$ alors il existe $t_{n,w} \in \bnrig{,r}{,K}$ tel que
$\iota_n(t_{n,w})=1 \mod{t^w K_n[\![t]\!]}$ et $\iota_m(t_{n,w}) \in t^w
K_n[\![t]\!]$ si $m \neq n$.
\end{lemm}

\begin{proof}
C'est une cons{\'e}quence imm{\'e}diate de la solution du probl{\`e}me des
{\og parties principales \fg} (voir \cite[\S8]{ML62}).
\end{proof}

Rappelons que les anneaux 
$\bnrig{,r}{,K}$ 
et donc aussi $\bnrig{}{,K}$ sont des
anneaux de B{\'e}zout, c'est-{\`a}-dire que tout id{\'e}al de type fini en
est principal. Ceci a un certain nombre de cons{\'e}quences pour
lesquelles on se reportera par exemple {\`a} \cite[\S 2]{KK00}.
Dans la suite, on pose $q=\phi(X)/X=((1+X)^p-1)/X$, ce qui fait que
$\iota_n(\phi^{n-1}(q))$ est une uniformisante de $K_n[\![t]\!]$.

\begin{prop}\label{idealdivt}
Si $I$ est un id{\'e}al principal de $\bnrig{,r}{,K}$, qui divise $(t^h)$
pour $h \geq 0$, alors
$I$ est engendr{\'e} par un {\'e}l{\'e}ment de la forme
$\prod_{n=n(r)}^{+\infty}(\phi^{n-1}(q)/p)^{j_n}$ avec $j_n \leq h$. 
\end{prop}

\begin{proof}
Rappelons que l'on a une d{\'e}composition $t = X \cdot 
\prod_{n=1}^{+\infty}(\phi^{n-1}(q)/p)$. Le lemme \cite[lemme
4.9]{LB02} montre que $\bnrig{,r}{,K}/ \phi^{n-1}(q) \simeq K_n$ et
donc que les id{\'e}aux $(\phi^{n-1}(q))$
sont premiers (et m{\^e}me maximaux) dans $\bnrig{,r}{,K}$.
Ceci montre que si $x$ divise $t^h$, 
alors $x$ est le produit d'une unit{\'e}
par un {\'e}l{\'e}ment de la forme
$\prod_{n=n(r)}^{+\infty}(\phi^{n-1}(q)/p)^{j_n}$ et il reste {\`a}
appliquer cela {\`a} un g{\'e}n{\'e}rateur de l'id{\'e}al $I$. 
\end{proof}

\subsection{Les $(\phi,\Gamma_K)$-modules sur l'anneau de
  Robba}\label{pgmodrap} 
Dans ce paragraphe, nous 
rappelons la d{\'e}finition des $\phi$-modules et des 
$(\phi,\Gamma_K)$-modules sur l'anneau de
Robba ainsi que quelques r{\'e}sultats
techniques concernant ces objets. 
Pour des d{\'e}finitions d'ordre plus g{\'e}n{\'e}ral, on peut voir
\cite[\S 2.5]{KK00}. Un  $\phi$-module sur $\bnrig{}{,K}$
est un $\bnrig{}{,K}$-module $\dbf$ libre de rang fini (not{\'e} $d$) et
muni d'une application $\phi$-semi-lin{\'e}aire toujours not{\'e}e 
$\phi: \dbf \ra \dbf$ telle que $\phi^*(\dbf)=\dbf$
(on note $\phi^*(\dbf)$ le $\phi$-module engendr{\'e} par $\phi(\dbf)$). Un
$(\phi,\Gamma_K)$-module sur $\bnrig{}{,K}$ est un $\phi$-module
muni en plus d'une action de $\Gamma_K$
semi-lin{\'e}aire par rapport l'action de ce groupe sur $\bnrig{}{,K}$ 
et commutant {\`a} $\phi$.  
On a un r{\'e}sultat de descente galoisienne (voir
aussi \cite[\S 2.5]{KK00}) pour les $(\phi,\Gamma_K)$-modules. 

\begin{defi}\label{actionglk}
Si $L/K$ est une extension galoisienne finie, et si $\dbf$ est un
$(\phi,\Gamma_L)$-module sur $\bnrig{}{,L}$, on dit que $\dbf$ est
muni d'une action de $G_{L/K}$ si le groupe $G_K$ agit sur $\dbf$ et
si de plus:
\begin{enumerate}
\item $H_L \subset G_K$ agit trivialement sur $\dbf$;
\item l'action de $G_L/H_L \subset G_K/H_L$ induite co{\"\i}ncide avec
celle de $\Gamma_L$.
\end{enumerate} 
\end{defi}

La proposition suivante se trouve essentiellement dans \cite[corollary
2.7]{KK00}:

\begin{prop}\label{descgal}
Si $L/K$ est une extension galoisienne finie et si $\dbf$ est un
$(\phi,\Gamma_L)$-module sur $\bnrig{}{,L}$ muni d'une action de
$G_{L/K}$, alors $\dbf^{H_K}$ est un $(\phi,\Gamma_K)$-module et 
$\dbf = \bnrig{}{,L} \otimes_{\bnrig{}{,K}} \dbf^{H_K}$.
\end{prop}

\begin{proof}
L'action de $\Gamma_K$ sur $\dbf^{H_K}$ est donn{\'e}e par
l'isomorphisme $\Gamma_K = G_K / H_K$. Le fait que
$\dbf = \bnrig{}{,L} \otimes_{\bnrig{}{,K}} \dbf^{H_K}$
suit de \cite[lemma 2.6]{KK00} (par exemple).
Montrons que $\phi^*(\dbf^{H_K})=(\dbf^{H_K})$. 
Si $x \in \dbf^{H_K}$, 
et si $\{y_i\}$ est une base de $\dbf$ sur $\bnrig{}{,L}$ 
contenue dans $\dbf^{H_K}$, alors on
peut {\'e}crire $x=\sum_{i=1}^d x_i \phi(y_i)$ avec $x_i \in
\bnrig{}{,L}$ et comme $x$ et les $y_i$ sont fixes par $H_K$, on a
aussi $x_i \in  \bnrig{}{,K}$ ce qui fait que $x \in
\phi^*(\dbf^{H_K})$. 
\end{proof}

Le th{\'e}or{\`e}me suivant est une variante d'un r{\'e}sultat de
Cherbonnier (voir \cite{FC96}). 

\begin{theo}\label{cherb}
Si $\dbf$ est un $\phi$-module sur $\bnrig{}{,K}$ alors il
existe $r(\dbf) \geq r(K)$ tel que pour tout $r \geq r(\dbf)$, il
existe un unique sous $\bnrig{,r}{,K}$-module $\dbf_r$ de $\dbf$
v{\'e}rifiant les propri{\'e}t{\'e}s suivantes:
\begin{enumerate}
\item $\dbf = \bnrig{}{,K} \otimes_{\bnrig{,r}{,K}} \dbf_r$;
\item le $\bnrig{,pr}{,K}$-module $\bnrig{,pr}{,K} \otimes_{\bnrig{,r}{,K}}
  \dbf_r$ a une base contenue dans $\phi(\dbf_r)$.
\end{enumerate}
En particulier, on a $\dbf_s= \bnrig{,s}{,K} \otimes_{\bnrig{,r}{,K}}
\dbf_r$ pour tous $s \geq r$ et si $\dbf$ est un
$(\phi,\Gamma_K)$-module, alors 
$\gamma(\dbf_r)=\dbf_r$ pour tout $\gamma \in \Gamma_K$. 
\end{theo}

\begin{proof}
Comme $\dbf$ est un $\bnrig{}{,K}$-module libre de rang $d$, il en
existe une base $e_1,\cdots,e_d$. Il existe alors $r=r(\dbf)$ tel que la
matrice de $\phi$ dans cette base est dans
$\on{GL}(d,\bnrig{,r}{,K})$ et si l'on pose $\dbf_r = \oplus_{i=1}^d
\bnrig{,r}{,K} e_i$, alors les deux conditions du th{\'e}or{\`e}me sont
remplies. Si $\dbf_r^{(1)}$ et $\dbf_r^{(2)}$ sont deux
$\bnrig{,r}{,K}$-modules satisfaisant les deux conditions ci-dessus,
et qu'on en choisit des bases, alors la matrice $M$ de
passage d'une base {\`a} l'autre et les matrices $P_1$ et $P_2$ de
$\phi$ dans ces deux bases vont satisfaire la relation
$\phi(M)=P_1^{-1}MP_2$ avec $P_1$ et $P_2 \in
\on{GL}(d,\bnrig{,pr}{,K})$ ce qui implique que $M \in
\on{M}(d,\bnrig{,r}{,K})$. En effet, si 
$M \in \on{M}(d,\bnrig{,s}{,K})$ avec $s > pr$, alors $P_1^{-1}MP_2$
et donc $\phi(M)$ est dans $\on{M}(d,\bnrig{,s}{,K})$. Mais
$\phi(\bnrig{,s}{,K}) \cap \bnrig{,s}{,K} = \phi(\bnrig{,s/p}{,K})$ ce
qui fait que si $s>pr$, alors on peut remplacer $s$ par $s/p$ et donc
finalement que $M \in \on{M}(d,\bnrig{,r}{,K})$.
Comme on peut en dire autant de $M^{-1}$,
c'est que $M \in \on{GL}(d,\bnrig{,r}{,K})$ et donc finalement que
$\dbf_r^{(1)} = \dbf_r^{(2)}$.
\end{proof}

{\'E}tant donn{\'e} un $\phi$-module $\dbf$ sur l'anneau de Robba, 
et $n \geq n(r)$ avec $r \geq r(\dbf)$, l'application $\iota_n :
\bnrig{,r}{,K} \ra K_n [\![t]\!]$ donne une structure de
$\iota_n(\bnrig{,r}{,K})$-module {\`a} $K_n [\![t]\!]$ et la formule
$\iota_n(\lambda) \cdot x = \lambda x$ donne une structure de
$\iota_n(\bnrig{,r}{,K})$-module {\`a} $\dbf_r$ 
que l'on note alors $\iota_n(\dbf_r)$, ce qui nous permet de
d{\'e}finir $K_n [\![t]\!] \otimes_{\iota_n(\bnrig{,r}{,K})}
\iota_n(\dbf_r)$. Pour all{\'e}ger les notations, on {\'e}crit plut{\^o}t: 
$K_n [\![t]\!] \otimes^{\iota_n}_{\bnrig{,r}{,K}} \dbf_r$.

\begin{prop}\label{unikiota}
Si $\dbf$ est un $\phi$-module de rang $d$
sur $\bnrig{}{,K}$ et si $\dbf^{(1)}$
et $\dbf^{(2)}$ sont deux sous-$\bnrig{}{,K}$-modules libres de rang
$d$ de $\dbf[1/t] = \bnrig{}{,K}[1/t] \otimes_{\bnrig{}{,K}} \dbf$
tels que:
\begin{enumerate}
\item $\dbf^{(i)}[1/t]= \dbf[1/t]$ si $i=1,2$;
\item $K_n [\![t]\!] \otimes^{\iota_n}_{\bnrig{,r}{,K}} \dbf^{(1)}_r  
= K_n [\![t]\!]
\otimes^{\iota_n}_{\bnrig{,r}{,K}} \dbf^{(2)}_r$ 
pour tout $n \gg 0$, 
\end{enumerate}
alors $\dbf^{(1)} = \dbf^{(2)}$. 
\end{prop}

\begin{proof}
Comme $\dbf^{(1)}[1/t] = \dbf^{(2)}[1/t]$, il existe $h \geq 0$ tel
que $t^h \dbf^{(2)} \subset \dbf^{(1)}$. Choisissons $r$ tel que $r
\geq \max(r(\dbf^{(1)}),r(\dbf^{(2)}))$ et tel qu'en plus  
\[ K_n [\![t]\!] \otimes^{\iota_n}_{\bnrig{,r}{,K}} \dbf^{(1)}_r  
= K_n [\![t]\!]
\otimes^{\iota_n}_{\bnrig{,r}{,K}} \dbf^{(2)}_r \] 
pour tout $n \geq n(r)$. Les diviseurs
{\'e}l{\'e}mentaires de $t^h \dbf^{(2)} \subset \dbf^{(1)}$ sont des
id{\'e}aux de $\bnrig{,r}{,K}$ qui divisent $t^h$ et donc de la
forme $(\prod_{n \geq n(r)} \phi^{n-1}(q^{\beta_{n,i}}))$ pour
$i=1,\hdots,d$ par la proposition \ref{idealdivt}. Comme le calcul des
diviseurs {\'e}l{\'e}mentaires commute {\`a} la localisation, on voit que
ceux de l'inclusion $t^h K_n [\![t]\!]
\otimes^{\iota_n}_{\bnrig{,r}{,K}} \dbf^{(2)}_r
\subset K_n [\![t]\!] \otimes^{\iota_n}_{\bnrig{,r}{,K}} \dbf^{(1)}_r$ 
sont donn{\'e}s par les id{\'e}aux $(t^{\beta_{n,i}})$ ce qui fait
que $\beta_{n,i}=h$ pour tous $n,i$ et donc finalement que $\dbf^{(1)} =
\dbf^{(2)}$. 
\end{proof}

\section{Construction de $(\phi,\Gamma_K)$-modules}\label{secII}
L'objet de ce chapitre est de montrer comment construire un
$\phi$-module sur  $\bnrig{}{,K}$, c'est-{\`a}-dire un objet
{\og global \fg}, {\`a} partir de conditions locales. Comme application,
on donne la construction d'un $(\phi,\Gamma_K)$-module sur
$\bnrig{}{,K}$ {\`a} partir d'un $(\phi,N,G_{L/K})$-module filtr{\'e}.

\subsection{Recollement de r{\'e}seaux locaux}\label{recres}
Si $\dbf$ est un $\phi$-module sur $\bnrig{}{,K}$ 
alors on a une application naturelle:
\[ \phi_n : K_{n+1}(\!(t)\!) \otimes_{K_n(\!(t)\!)} 
\left[ K_n(\!(t)\!) \otimes^{\iota_n}_{\bnrig{,r}{,K}}
\dbf_r \right]  \lra K_{n+1}(\!(t)\!) 
\otimes^{\iota_{n+1}}_{\bnrig{,r}{,K}} \dbf_r \]
d{\'e}finie par $\phi_n [ f(t) \otimes (g(t) 
\otimes \iota_n(x))] = f(t)g(t) \otimes
\iota_{n+1} (\phi(x))$ en utilisant le fait que $\phi(\dbf_r) \subset
\bnrig{,pr}{,K} \otimes_{\bnrig{,r}{,K}} \dbf_r$, que
$\iota_{n+1}(\bnrig{,pr}{,K}) \subset K_{n+1}[\![t]\!]$ et que
$\iota_{n+1}(\phi(x)) = \iota_n(x)$. 

\begin{defi}\label{deficomp}
Si $\dbf$ est un $\phi$-module sur $\bnrig{}{,K}$ et si 
$r \geq r(\dbf)$ et 
$\{M_n\}_{n \geq n(r)}$ est une suite de 
$K_n[\![t]\!]$-r{\'e}seaux de $K_n(\!(t)\!) 
\otimes^{\iota_n}_{\bnrig{,r}{,K}} \dbf_r$, alors on dit 
que $\{M_n\}$ est une suite $\phi$-compatible
si $\phi_n(K_{n+1}[\![t]\!] \otimes_{K_n[\![t]\!]} M_n) = M_{n+1}$.  
\end{defi}

On v{\'e}rifie sans mal que la donn{\'e}e d'un
sous-$\phi$-module $\mbf$ 
de rang maximal d'un $\phi$-module
$\dbf$ d{\'e}termine une suite $\phi$-compatible de 
$K_n[\![t]\!]$-r{\'e}seaux de $K_n(\!(t)\!) 
\otimes^{\iota_n}_{\bnrig{,r}{,K}} \dbf_r$ en posant 
$M_n= K_n[\![t]\!] \otimes^{\iota_n}_{\bnrig{,r}{,K}} \mbf_r$. Le
th{\'e}or{\`e}me suivant donne la r{\'e}ciproque de cette construction.

\begin{theo}\label{constm}
Si $\dbf$ est un $\phi$-module sur $\bnrig{}{,K}$ et si
$\{M_n\}_{n \geq n(r)}$ est une suite 
$\phi$-compatible de r{\'e}seaux de $K_n(\!(t)\!)
\otimes^{\iota_n}_{\bnrig{,r}{,K}} \dbf_r$, 
alors il existe un unique sous
$\phi$-module $\mbf$ de $\dbf[1/t]$ tel que $K_n[\![t]\!]
\otimes^{\iota_n}_{\bnrig{,r}{,K}} \mbf_r = M_n$ pour 
tout $n \geq n(r)$. Enfin, on a $\mbf[1/t]=\dbf[1/t]$. 

Si $\dbf$ est un $(\phi,\Gamma_K)$-module et si les $M_n$ sont
stables sous l'action induite de $\Gamma_K$, alors $\mbf$ est lui
aussi un $(\phi,\Gamma_K)$-module.
\end{theo}

L'unicit{\'e} de $\mbf$ r{\'e}sulte imm{\'e}diatement de la proposition
\ref{unikiota} et le reste de ce paragraphe est consacr{\'e} {\`a} la
d{\'e}monstration de l'existence d'un tel $\mbf$.

\begin{lemm}\label{res1}
Il existe un entier $h \geq 0$ tel que pour tout $n \geq n(r)$ on ait 
\[ t^h K_n[\![t]\!] \otimes^{\iota_n}_{\bnrig{,r}{,K}} \dbf_r
\subset M_n \subset t^{-h} K_n[\![t]\!]
\otimes^{\iota_n}_{\bnrig{,r}{,K}} \dbf_r. \]
\end{lemm}

\begin{proof}
Comme $M_{n(r)}$ est un r{\'e}seau de $K_{n(r)}(\!(t)\!) 
\otimes^{\iota_{n(r)}}_{\bnrig{,r}{,K}}
\dbf_r$, il existe $h$ tel que l'inclusion ci-dessus est
vraie pour $n=n(r)$. Le fait que \[ \phi_n(K_{n+1}[\![t]\!]
\otimes_{K_n[\![t]\!]} K_n[\![t]\!] \otimes^{\iota_n}_{\bnrig{,r}{,K}}
\dbf_r) = K_{n+1}[\![t]\!] \otimes^{\iota_{n+1}}_{\bnrig{,r}{,K}}
\dbf_r \] 
et que $\phi_n(K_{n+1}[\![t]\!] \otimes_{K_n[\![t]\!]} M_n)=M_{n+1}$
montrent que si l'inclusion est vraie pour $n$ elle est aussi vraie
pour $n+1$, ce qui montre le lemme par r{\'e}currence.
\end{proof}

\begin{lemm}\label{res2}
Si on pose $\mbf_r = \{ x \in
t^{-h} \dbf_r$ tels que $\iota_n(x) \in M_n$ pour tout $n \geq n(r) \}$,
alors $\mbf_r$ est un $\bnrig{,r}{,K}$-module libre de rang $d$. 
\end{lemm}

\begin{proof}
Comme les applications $\iota_n : t^{-h} \dbf_r \ra M_n[1/t]$ sont
continues, $\mbf_r$ est ferm{\'e} dans $t^{-h} \dbf_r$. D'autre part, le
lemme \ref{res1} montre que $t^h \dbf_r \subset \mbf_r$ et le
th{\'e}or{\`e}me de Forster (cf \cite[th{\'e}or{\`e}me 4.10]{LB02} pour une
d{\'e}monstration) montre alors que $\mbf_r$ est un
$\bnrig{,r}{,K}$-module libre de rang $d$. 
\end{proof}

\begin{lemm}\label{res3}
On a $K_n[\![t]\!] \otimes^{\iota_n}_{\bnrig{,r}{,K}} \mbf_r = M_n$ 
pour tout $n \geq n(r)$.   
\end{lemm}

\begin{proof}
Comme $K_n[\![t]\!] \otimes^{\iota_n}_{\bnrig{,r}{,K}} \mbf_r$ 
et $M_n$ sont complets
pour la topologie $t$-adique, il suffit de montrer que l'application
naturelle $K_n[\![t]\!] \otimes^{\iota_n}_{\bnrig{,r}{,K}} 
\mbf_r \ra M_n / t^h M_n$ est
surjective pour tout $n$.    
Si $x \in M_n$, alors le lemme \ref{res1} montre qu'il existe $y \in
t^{-h} \dbf_r$ tel que $\iota_n(y) - x \in t^h M_n$. 
Le lemme \ref{partunit} appliqu{\'e} {\`a} $w=3h$ nous donne un
{\'e}l{\'e}ment $t_{n,3h}$ et on pose
$z = t_{n,3h} y$ ce qui fait que
\[ \iota_n(z) - \iota_n(y) \in t^{2h} 
K_n[\![t]\!] \otimes^{\iota_n}_{\bnrig{,r}{,K}} \dbf_r
\subset t^h M_n \] 
tandis que si $m \neq n$ alors 
\[ \iota_m(z) \in t^{2h} 
K_m[\![t]\!] \otimes^{\iota_m}_{\bnrig{,r}{,K}} \dbf_r 
\subset t^h M_m \subset M_m \]
ce qui fait 
finalement que $z \in \mbf_r$ et on trouve bien que l'application
naturelle \[ K_n[\![t]\!] \otimes^{\iota_n}_{\bnrig{,r}{,K}} 
\mbf_r \ra M_n / t^h M_n \] 
est surjective.     
\end{proof}

\begin{proof}[D{\'e}monstration du th{\'e}or{\`e}me \ref{constm}]
Les deux lemmes pr{\'e}c{\'e}dents montrent que si l'on pose  
$\mbf_r = \{ x \in
t^{-h} \dbf_r$ tels que $\iota_n(x) \in M_n$ pour tout $n \geq n(r) \}$,
alors $\mbf_r$ est un $\bnrig{,r}{,K}$-module libre de rang $d$ et
$K_n[\![t]\!] \otimes^{\iota_n}_{\bnrig{,r}{,K}} 
\mbf_r = M_n$ pour tout $n \geq
n(r)$ et qu'on peut donc poser
$\mbf = \bnrig{}{,K} \otimes_{\bnrig{,r}{,K}} \mbf_r$. 
Le fait que la suite de r{\'e}seaux $\{M_n\}$ est 
$\phi$-compatible montre que
$\phi^*(\mbf) \subset \mbf$ et que $K_n[\![t]\!] 
\otimes^{\iota_n}_{\bnrig{,r}{,K}} 
\phi^*(\mbf)_{pr} = M_n$ pour tout $n \geq
n(pr)$ et la proposition \ref{unikiota} 
appliqu{\'e}e {\`a} $\dbf^{(1)}=\mbf$ et $\dbf^{(2)}=\phi^*(\mbf)$
montre alors qu'en fait $\phi^*(\mbf)=\mbf$. Enfin le fait que 
$\mbf[1/t]=\dbf[1/t]$ suit imm{\'e}diatement du lemme \ref{res2}.

L'assertion quant {\`a} l'action {\'e}ventuelle de $\Gamma_K$ est
{\'e}vidente. 
\end{proof}

\subsection{Des $(\phi,N,G_{L/K})$-modules filtr{\'e}s aux
  $(\phi,\Gamma_K)$-modules}\label{pntopgm}
Soit $\ell_X$ une variable; 
on prolonge les actions de $\phi$ et de
$\Gamma_K$ {\`a} $\bnrig{}{,K}[\ell_X]$ par les formules suivantes:
$\phi(\ell_X) = p \ell_X + \log(\phi(X)/X^p)$ et
$\gamma(\ell_X)=\ell_X+\log(\gamma(X)/X)$. Bien s{\^u}r, il faut penser
{\`a} $\ell_X$ comme {\`a} {\og $\log(X)$ \fg}. On d{\'e}finit alors un
op{\'e}rateur de monodromie $N$ sur $\bnrig{}{,K}[\ell_X]$ par la
formule $N(\ell_X)=-p/(p-1)$ et on prolonge l'application $\iota_n$
{\`a} $\bnrig{,r}{,K}[\ell_X]$ par
$\iota_n(\ell_X)=\log(\eps^{(n)}\exp(t/p^n)-1) \in K_n[\![t]\!]$.

Si $D$ est un $(\phi,N)$-module filtr{\'e} (relatif {\`a} $K$), 
on pose $\dbf = (\bnrig{}{,K}[\ell_X] \otimes_F D)^{N=0}$.  
Pour tout $n \in \ZZ$, on a $K_0=\phi^{-n}(K_0) \subset K$ 
ce qui nous donne une structure de $\phi^{-n}(K_0)$-module sur $K$ et
sur $D$ que l'on note alors $\iota_n(D)$ et nous
{\'e}crivons $K \otimes^{\iota_n}_{K_0} D$ au lieu de  
$K \otimes_{\phi^{-n}(K_0)} \iota_n(D)$ pour all{\'e}ger les notations;
l'application $\xi_n : K
\otimes_{K_0} D \ra K \otimes^{\iota_n}_{K_0} D$ 
qui envoie $\mu \otimes x$ sur $\mu \otimes \iota_n(\phi^n(x))$  
est un alors isomorphisme (rappelons que $\iota_n=\phi^{-n}$) 
que l'on utilise pour d{\'e}finir une filtration sur $D^n_K = K
\otimes^{\iota_n}_{K_0} D$. 
On d{\'e}finit une filtration sur $K_n(\!(t)\!)$ par la formule
$\fil^i K_n(\!(t)\!) = t^i K_n[\![t]\!]$ ce qui nous donne une
filtration sur $K_n(\!(t)\!) \otimes_K D^n_K$, et on pose $M_n(D) =
\fil^0 (K_n(\!(t)\!) \otimes_K D^n_K)$. Le foncteur $D \mapsto M_n(D)$
est alors un $\otimes$-foncteur exact. 

\begin{prop}
La famille de r{\'e}seaux $\{M_n\}$ de $K_n(\!(t)\!)
\otimes^{\iota_n}_{\bnrig{,r}{,K}} \dbf_r$ d{\'e}finie par 
$M_n(D) = \fil^0 (K_n(\!(t)\!) \otimes_K D^n_K)$ 
est $\phi$-compatible. 
\end{prop}

\begin{proof}
Le $K_n[\![t]\!]$-module $M_n(D) = \fil^0
(K_n(\!(t)\!) \otimes_K D_K^n)$ est libre de rang $d$, 
engendr{\'e} par une base
de la forme $t^{-h_i} \otimes \xi_n(e_i)$ 
o{\`u} $\{e_i\}_{1 \leq i \leq d}$ est une
base de $D_K$ adapt{\'e}e {\`a} la filtration et $h_i = t_H(e_i)$, ce qui
fait que $M_{n+1}(D) = \phi_n(K_{n+1}[\![t]\!] \otimes_{K_n[\![t]\!]}
M_n(D))$ puisque $\xi_{n+1} = \phi_n \circ \xi_n$ sur $D_K$. 
\end{proof}

\begin{defi}
Si $D$ est un $(\phi,N)$-module filtr{\'e} relatif {\`a} $K$, soit $\pgm(D)$ le
$(\phi,\Gamma_K)$-module fourni par le th{\'e}or{\`e}me \ref{constm} {\`a} partir
des r{\'e}seaux $M_n(D) = \fil^0 (K_n(\!(t)\!) \otimes_K D_K^n)$ de $\dbf = 
(\bnrig{}{,K}[\ell_X] \otimes_F D)^{N=0}$. 
\end{defi}

\begin{prop}
Si $D$ est un $(\phi,N,G_{L/K})$-module filtr{\'e} et $D'$ est le
$(\phi,N)$-module filtr{\'e} relatif {\`a} $L$ 
qu'on en d{\'e}duit (par oubli de l'action de $G_{L/K}$), 
alors $\pgm(D')$ est
un $(\phi,\Gamma_L)$-module sur $\bnrig{,r}{,L}$ muni d'une action de
$G_{L/K}$ (cf d{\'e}finition \ref{actionglk}).
\end{prop}

\begin{proof}
V{\'e}rification imm{\'e}diate.
\end{proof}

\begin{defi}
Si $D$ est un $(\phi,N,G_{L/K})$-module filtr{\'e}, alors on d{\'e}finit
$\pgm(D) = \pgm(D')^{H_K}$.
\end{defi}

Rappelons (cf \cite[4.3.2]{Bu88sst}) que par
d{\'e}finition on a une filtration $\fil^i D_L$ sur $D_L = L
\otimes_{L_0} D$ et que si $\fil^i D_K$ est la filtration 
induite sur $D_K = K \otimes_{K_0} D$, alors $\fil^i D_L = L \otimes_K 
\fil^i D_K$. On a construit au d{\'e}but du paragraphe un isomorphisme
$G_{L/K}$-{\'e}quivariant $D_L \ra D_L^n$ que l'on utilise pour
d{\'e}finir $D_K^n = (D_L^n)^{G_{L/K}}$ ce qui fait que $D_L^n = L
\otimes_K D_K^n$.

\begin{prop}\label{mndk}
Si $D$ est un $(\phi,N,G_{L/K})$-module filtr{\'e}, et si on pose
$M_n(D_K)=\fil^0 (K_n(\!(t)\!) \otimes_K D_K^n)$, alors
$K_n[\![t]\!] \otimes^{\iota_n}_{\bnrig{,r}{,K}} 
\pgm(D)_r = M_n(D_K)$ pour tout $n \geq n(r)$.
\end{prop}

\begin{proof}
Par construction, $K_n[\![t]\!] 
\otimes^{\iota_n}_{\bnrig{,r}{,K}} \pgm(D)_r \subset 
\fil^0 (L_n(\!(t)\!) \otimes_L D_L^n)^{H_K}$. 
Ce que l'on a rappel{\'e}
ci-dessus montre que $\fil^0 (L_n(\!(t)\!) \otimes_L D_L^n)^{H_K} = \fil^0
(K_n(\!(t)\!) \otimes_K D_K^n)$ et donc 
$K_n[\![t]\!] \otimes^{\iota_n}_{\bnrig{,r}{,K}}
\pgm(D)_r \subset M_n(D_K)$. 

Si $x \in M_n(D_K)$, alors $x \in L_n[\![t]\!] 
\otimes^{\iota_n}_{\bnrig{,r}{,L}}
\pgm(D')_r  = L_n[\![t]\!] \otimes^{\iota_n}_{\bnrig{,r}{,K}}
\pgm(D)_r$. Finalement, comme $x$ est fix{\'e} par $H_K$, on a 
$x \in (L_n[\![t]\!] \otimes^{\iota_n}_{\bnrig{,r}{,K}}
\pgm(D)_r)^{H_K} = K_n[\![t]\!] \otimes^{\iota_n}_{\bnrig{,r}{,K}}
\pgm(D)_r$ ce qui montre que l'application 
$K_n[\![t]\!] \otimes^{\iota_n}_{\bnrig{,r}{,K}}
\pgm(D)_r \ra M_n(D_K)$ est un isomorphisme.
\end{proof}

\begin{theo}\label{eqcat}
Le foncteur $D \mapsto \pgm(D)$ est un $\otimes$-foncteur exact 
de la cat{\'e}gorie des $(\phi,N,G_K)$-modules
filtr{\'e}s dans la cat{\'e}gorie des $(\phi,\Gamma_K)$-modules sur
$\bnrig{}{,K}$ et le rang de $\pgm(D)$ est {\'e}gal {\`a} la dimension de
$D$. 
\end{theo}

\begin{proof}
Si l'on a une suite exacte $0 \ra D_1 \ra D_2
\ra D_3 \ra 0$, alors montrons que
$\pgm(D_2) \ra \pgm(D_3)$ est surjectif. Comme on l'a dit plus haut,  
les foncteurs $D \mapsto
M_n(D)$ sont des $\otimes$-foncteurs exacts.
Comme $M_n(D_2) \ra M_n(D_3)$ est surjectif pour tout $n$,
l'image $\mbf$ 
de $\pgm(D_2)$ dans $\pgm(D_3)$ v{\'e}rifie la condition que
$K_n[\![t]\!] \otimes^{\iota_n}_{\bnrig{,r}{,K}} \mbf_r = M_n(D_3)$ 
du th{\'e}or{\`e}me
\ref{constm} pour $D=D_3$ et par unicit{\'e}, 
on a donc que l'image de 
$\pgm(D_2)$ dans $\pgm(D_3)$ est $\pgm(D_3)$ tout entier, ce qui fait
que $\pgm(D_2) \ra \pgm(D_3)$ est bien surjectif.
La v{\'e}rification du fait que $\pgm(D_1
\otimes D_2)=\pgm(D_1) \otimes \pgm(D_2)$ et celles des autres
conditions sont similaires. 
\end{proof}

Dans le chapitre suivant, nous allons d{\'e}terminer l'image essentielle
du foncteur $D \mapsto \pgm(D)$ et montrer que ce foncteur est une
{\'e}quivalence de cat{\'e}gories sur son image essentielle.

\section{Construction de $(\phi,N,G_K)$-modules filtr{\'e}s}\label{chapmonod}
L'objet de ce chapitre est de montrer comment on peut associer {\`a}
certains $(\phi,\Gamma_K)$-modules sur $\bnrig{}{,K}$ un
$(\phi,N,G_K)$-module filtr{\'e}. Cette construction est un inverse de
celle du chapitre pr{\'e}c{\'e}dent.

\subsection{L'alg{\`e}bre de Lie de $\Gamma_K$}\label{liealg}
Le groupe $\Gamma_K$ s'identifie, via le caract{\`e}re cyclotomique, 
{\`a} un sous-groupe ouvert de $\Zp^*$, et c'est donc un groupe de Lie
$p$-adique de dimension $1$. Son alg{\`e}bre de Lie: $\on{Lie}(\Gamma_K)$
est un $\Qp$-espace vectoriel de rang $1$ dont une base est donn{\'e}e
par l'op{\'e}rateur $\log(\gamma)/\log_p \chi(\gamma)$ qui ne d{\'e}pend
pas du choix de $\gamma \in \Gamma_K$. 

\begin{prop}\label{deficonn}
Si $\dbf$ est un $(\phi,\Gamma_K)$-module sur $\bnrig{}{,K}$ alors
l'alg{\`e}bre de Lie de $\Gamma_K$ agit par un
op{\'e}rateur diff{\'e}rentiel $\nabla_{\dbf} : \dbf \ra \dbf$ qui commute {\`a}
$\phi$ et {\`a} l'action de $\Gamma_K$ et qui v{\'e}rifie $\nabla_{\dbf} (\lambda
x)  = \nabla(\lambda) x + \lambda \nabla_{\dbf}(x)$.
\end{prop}

\begin{proof}
La d{\'e}monstration est semblable {\`a} celle de \cite[lemme 5.2]{LB02}
et du paragraphe qui la suit.
Soit $V_{[r;s]}$ la valuation sup sur la couronne $C_K[r;s]$. La topologie de
$\dbf_r$ est la topologie de Fr{\'e}chet d{\'e}finie par l'ensemble
$\{V_{[r;s]}\}_{s \geq r}$. Fixons $s \geq r$; l'action de $\Gamma_K$
sur $\dbf_r$ est continue, et il existe donc 
un sous-groupe ouvert $\Gamma_s \subset \Gamma_K$ tel
que $V_{[r;s]}((1-\gamma)x) \geq V_{[r;s]}(x)+1$ pour tout $x \in
\dbf_r$ et $\gamma \in \Gamma_s$. La s{\'e}rie d'op{\'e}rateurs
\[  - \frac{1}{\log_p \chi(\gamma)} \sum_{n=1}^{\infty}
\frac{(1-\gamma)^n}{n} \] converge alors vers un op{\'e}rateur continu
$\nabla_{\dbf,s}$ de $\dbf_r$ vers sa compl{\'e}tion pour $V_{[r;s]}$,
qui ne d{\'e}pend pas du choix de $\gamma \in \Gamma_s$. Ces
op{\'e}rateurs $\nabla_{\dbf,s}$ se recollent en un op{\'e}rateur
$\nabla_{\dbf} :
\dbf_r \ra \dbf_r$ qui est continu pour la topologie de Fr{\'e}chet. 

Enfin le fait que $\nabla_{\dbf}(\lambda
x)  = \nabla(\lambda) x + \lambda \nabla_{\dbf}(x)$ r{\'e}sulte par passage
{\`a} la limite du fait que
\[ (1-\gamma)(\lambda \cdot x) =
(1-\gamma)(\lambda) \cdot x+\lambda \cdot (1-\gamma)(x) - 
(1-\gamma)(\lambda) \cdot (1-\gamma)(x). \]
\end{proof}

\begin{defi}\label{loctriv}
On dit que l'op{\'e}rateur $\nabla_{\dbf}$ est localement trivial sur un
$(\phi,\Gamma_K)$-module $\dbf$ s'il existe $r$ tel que 
\[ K_{n(r)}(\!(t)\!) \otimes^{\iota_{n(r)}}_{\bnrig{,r}{,K}} \dbf_r 
= K_{n(r)}(\!(t)\!) \otimes_{K_{n(r)}} \left( K_{n(r)}(\!(t)\!) 
\otimes^{\iota_{n(r)}}_{\bnrig{,r}{,K}} \dbf_r
\right)^{\nabla_{\dbf}=0}.  \]
\end{defi}

Avant de continuer, faisons quelques rappels sur les modules {\`a}
connexion. Soit $E$ un corps de caract{\'e}ristique $0$,
et $M$ un $E(\!(t)\!)$-espace vectoriel 
de dimension finie $d$ muni
d'une connexion $\nabla_M : M \ra M$ qui {\'e}tend $\nabla : f(t)
\mapsto t \cdot df/dt$. 
Une section horizontale de $M$ est un {\'e}l{\'e}ment de $M^{\nabla_M=0}$
et un argument classique montre que $\dim_E M^{\nabla_M=0} \leq d$.

On dit que la connexion $\nabla_M$ est
r{\'e}guli{\`e}re si $M$ poss{\`e}de un $E[\![t]\!]$-r{\'e}seau $M_0$ tel que
$\nabla_M(M_0) \subset M_0$ et on dit que la connexion $\nabla_M$ est
triviale si $M$ poss{\`e}de un $E[\![t]\!]$-r{\'e}seau $M_0$ tel que
$\nabla_M(M_0) \subset t M_0$. Un argument d'approximations successives
montre que dans ce cas, $M_0^{\nabla_M=0}$ est un $E$-espace vectoriel
de dimension $d$ et que $M_0 = E[\![t]\!] \otimes_E M_0^{\nabla_M=0}$. On
voit donc que la connexion $\nabla_M$ est triviale si et seulement si
$\dim_E M^{\nabla_M=0} = d$, et que dans ce cas $M_0 = E[\![t]\!] \otimes_E
M^{\nabla_M=0}$ est l'unique $E[\![t]\!]$-r{\'e}seau de $M$ tel que
$\nabla_M(M_0) \subset t M_0$. 

\begin{lemm}\label{soustriv}
Si $N$ est un sous-$E(\!(t)\!)$-espace vectoriel de $M$ stable par la
connexion $\nabla_M$, et si $\nabla_M$ est triviale sur $M$, 
alors elle est aussi triviale sur $N$.
\end{lemm}

\begin{proof}
Si $M_0$ est un $E[\![t]\!]$-r{\'e}seau de $M$ tel que
$\nabla_M(M_0) \subset t M_0$, alors $M_0 \cap N$ est un
$E[\![t]\!]$-r{\'e}seau de $N$ et $\nabla_M(M_0 \cap N)
\subset t (M_0 \cap N)$ ce qui fait que la connexion 
$\nabla_M$ est triviale sur $N$. 
\end{proof}

La terminologie de la d{\'e}finition \ref{loctriv} est compatible
avec ce que l'on vient de rappeler: $\nabla_{\dbf}$ est localement
triviale si et seulement si $\nabla_{\dbf}$ est triviale sur le
$K_{n(r)}(\!(t)\!)$-module {\`a} connexion $K_{n(r)}(\!(t)\!)
\otimes^{\iota_{n(r)}}_{\bnrig{,r}{,K}} \dbf_r$. En particulier, si 
$\nabla_{\dbf}$ est localement triviale sur $\dbf$ et si $\dbf'$ est
un sous-objet de $\dbf$, alors $\nabla_{\dbf}$ est localement triviale
sur $\dbf'$. 

\begin{lemm}\label{trivdivers}
Si $\nabla_{\dbf}$ est localement triviale sur $\dbf$, 
alors, $\nabla_{\dbf}$ est triviale sur $K_n(\!(t)\!)
\otimes^{\iota_n}_{\bnrig{,r}{,K}} \dbf_r$ pour tout $n \geq n(r)$. 

Dans ce cas, si $D_n$ est le r{\'e}seau de $K_n(\!(t)\!)
\otimes^{\iota_n}_{\bnrig{,r}{,K}} \dbf_r$ tel que $\nabla_{\dbf}(D_n)
\subset t D_n$, alors $D_{n+1}=\phi_n (K_{n+1}[\![t]\!] \otimes_{K_n[\![t]\!]}
D_n)$. 
\end{lemm}

\begin{proof}
Si $D_n =  K_n[\![t]\!] \otimes_{K_n} ( K_n(\!(t)\!) 
\otimes^{\iota_n}_{\bnrig{,r}{,K}} \dbf_r )^{\nabla_{\dbf}=0}$,
alors par hypoth{\`e}se, $\nabla_{\dbf}(D_{n(r)}) \subset t D_{n(r)}$ et
si $n \geq n(r)$ et $\nabla_{\dbf}(D_n) \subset t D_n$ alors
$\nabla_{\dbf} (D_{n+1}) \subset t D_{n+1}$ si 
$D_{n+1}$ est le r{\'e}seau de $K_{n+1}(\!(t)\!)
\otimes^{\iota_{n+1}}_{\bnrig{,r}{,K}} \dbf_r$ donn{\'e} par 
$D_{n+1}=\phi_n (K_{n+1}[\![t]\!] \otimes_{K_n[\![t]\!]} D_n)$ 
puisque $\nabla_{\dbf}$ commute {\`a} $\phi_n$, 
ce qui montre le r{\'e}sultat par r{\'e}currence. 
\end{proof}

\begin{prop}\label{mbfistriv}
Si $D$ est un $(\phi,N,G_{L/K})$-module filtr{\'e}, alors $\nabla_{\dbf}$ est
localement triviale sur $\pgm(D)$.
\end{prop}

\begin{proof}
Par la proposition \ref{mndk}, on a 
$K_n(\!(t)\!) \otimes^{\iota_n}_{\bnrig{,r}{,K}} 
\pgm(D)_r = K_n(\!(t)\!) \otimes_K D_K^n$ et comme $G_{L/K}$ est fini,
$\nabla_{\dbf}=0$ sur $D_K^n$.
\end{proof}

Cette proposition montre que si $\mbf$ est dans l'image du foncteur $D
\mapsto \pgm(D)$, alors $\nabla_{\mbf}$ est
localement triviale sur $\mbf$. Nous allons voir que r{\'e}ciproquement,
si $\nabla_{\mbf}$ est
localement triviale sur $\mbf$, alors $\mbf$ est dans l'image 
du foncteur $D \mapsto \pgm(D)$.

\subsection{{\'E}quations diff{\'e}rentielles $p$-adiques}\label{eqdiff}
Dans ce paragraphe, nous d{\'e}terminons l'image essentielle du foncteur
$D \mapsto \pgm(D)$. Si $\dbf$ est un $\phi$-module sur $\bnrig{}{,K}$ muni d'un
op{\'e}rateur diff{\'e}rentiel $\partial_{\dbf}$ qui {\'e}tend l'op{\'e}rateur
$\partial: f \mapsto (1+X)df/dX$ et tel que 
$\partial_{\dbf} \circ \phi = p \cdot \phi \circ \partial_{\dbf}$, 
alors on  dit que $\dbf$ est une
{\'e}quation diff{\'e}rentielle $p$-adique avec structure de Frobenius.

Si $\dbf$ est un $(\phi,\Gamma_K)$-module 
sur $\bnrig{}{,K}$ et si $\dbf$ est stable par l'op{\'e}rateur 
$\partial_{\dbf} = t^{-1} \nabla_{\dbf}$, alors $\dbf$ est une
{\'e}quation diff{\'e}rentielle $p$-adique avec structure de Frobenius.

Rappelons tout d'abord le th{\'e}or{\`e}me de monodromie $p$-adique,
conjectur{\'e} par Crew et
d{\'e}montr{\'e} par Andr{\'e}, Kedlaya et Mebkhout 
(cf \cite{YA01}, \cite{KK00} et \cite{ZM01} ainsi que 
le s{\'e}minaire Bourbaki \cite{PC01}):

\begin{theo}\label{akm}
Si $\dbf$ est une {\'e}quation diff{\'e}rentielle $p$-adique avec
structure de Frobenius, 
alors il existe une extension finie $L/K$ telle que l'application
naturelle 
\[ \bnrig{}{,L}[\ell_X] \otimes_{L_0'} \left( \bnrig{}{,L}[\ell_X]
\otimes_{\bnrig{}{,K}} \dbf \right)^{\partial_D = 0} \ra 
\bnrig{}{,L}[\ell_X] \otimes_{\bnrig{}{,K}} \dbf \]
est un isomorphisme.
\end{theo}

Supposons maintenant que $\dbf$ est un $(\phi,\Gamma_K)$-module 
sur $\bnrig{}{,K}$ qui est stable par l'op{\'e}rateur 
$\partial_{\dbf} = t^{-1} \nabla_{\dbf}$ et
posons $S_L(\dbf)=(\bnrig{}{,L}[\ell_X] 
\otimes_{\bnrig{}{,K}} \dbf)^{\partial_D = 0}$. C'est un $L_0'$-espace
vectoriel qui h{\'e}rite d'une action r{\'e}siduelle de $\Gamma_K$
triviale sur un sous-groupe ouvert (puisque $t \partial_D =
0$). Quitte {\`a} remplacer $L$ par une extension finie, on peut donc
supposer que $\Gamma_L$ agit trivialement sur $S_L(\dbf)$, ce qui fait
que $L_0'=L_0$ et on pose alors $\on{Sol}_L(\dbf)=S_L(\dbf)$
%
%
ce qui fait que \[ \on{Sol}_L(\dbf) = \left(\bnrig{}{,L}[\ell_X] 
\otimes_{\bnrig{}{,K}} \dbf \right)^{G_L}\] 
est un $(\phi,N,G_{L/K})$-module tel que 
\[ \bnrig{}{,L}[\ell_X] \otimes_{L_0}  \on{Sol}_L(\dbf) =  
\bnrig{}{,L}[\ell_X] \otimes_{\bnrig{}{,K}} \dbf. \]

\begin{defi}\label{defisol}
Le $L_0$-espace vectoriel $\on{Sol}_L(\dbf)$ est alors appel{\'e} l'espace des
$G_L$-solutions du $(\phi,\Gamma_K)$-module $\dbf$.
\end{defi}

On peut donc reformuler le th{\'e}or{\`e}me \ref{akm} ci-dessus et la
discussion qui suit en disant que si $\dbf$ est un 
$(\phi,\Gamma_K)$-module sur $\bnrig{}{,K}$ 
qui est stable par l'op{\'e}rateur 
$\partial_{\dbf} = t^{-1} \nabla_{\dbf}$, 
alors il admet des $G_L$-solutions pour
$L$ assez grand.

\begin{theo}\label{imageessent}
Si $\mbf$ est un $(\phi,\Gamma_K)$-module sur $\bnrig{}{,K}$ tel que
$\nabla_{\mbf}$ est localement triviale, alors il existe un
unique $(\phi,\Gamma_K)$-module $\dbf \subset \mbf[1/t]$ tel que
$\dbf[1/t]=\mbf[1/t]$ et tel que $\partial_{\mbf} (\dbf) \subset \dbf$. 

De plus, la donn{\'e}e de $\mbf$ d{\'e}termine une filtration sur $L
\otimes_{L_0} \on{Sol}_L(\dbf)$ et donc une structure de
$(\phi,N,G_{L/K})$-module filtr{\'e} sur $\on{Sol}_L(\dbf)$
qui a la propri{\'e}t{\'e} que $\mbf = \pgm(\on{Sol}_L(\dbf))$. 
\end{theo}

\begin{proof}
Comme $\nabla_{\mbf}$ est localement triviale, il existe une famille de
r{\'e}seaux $D_n$ de $K_n(\!(t)\!) \otimes^{\iota_n}_{\bnrig{,r}{,K}}
\mbf_r$ tels que $\nabla_{\dbf}(D_n) \subset t D_n$ et  
comme $D_{n+1} = \phi_n (K_{n+1}[\![t]\!] \otimes_{K_n[\![t]\!]} D_n)$ 
par le lemme \ref{trivdivers}, la famille $\{D_n\}$ est
$\phi$-compatible. Par le th{\'e}or{\`e}me \ref{constm}, il existe donc un
$(\phi,\Gamma_K)$-module $\dbf$ tel que $K_n[\![t]\!]
\otimes^{\iota_n}_{\bnrig{,r}{,K}} 
\dbf_r = D_n$. Comme $\nabla_{\dbf}(D_n) \subset t D_n$ pour tout
$n$, on a $\nabla_{\dbf}(\dbf) \subset t \dbf$ ce qui fait que $\dbf$ muni de
la connexion $\partial_{\dbf} = t^{-1} \nabla_{\dbf}$ est une {\'e}quation
diff{\'e}rentielle $p$-adique avec structure de Frobenius
ce qui montre le premier point. L'unicit{\'e} de $\dbf$ suit du fait que
l'on a n{\'e}cessairement $K_n[\![t]\!] \otimes^{\iota_n}_{\bnrig{,r}{,K}}
\dbf_r = D_n$ pour tout $n \geq n(r)$. 

Par le th{\'e}or{\`e}me de monodromie $p$-adique de Andr{\'e}, Kedlaya et
Mebkhout que l'on a rappel{\'e} ci-dessus en \ref{akm}, 
et la discussion qui le suit, 
le $(\phi,\Gamma_K)$-module $\dbf$
admet des $G_L$-solutions pour $L/K$ assez grand. 
Posons $D=\on{Sol}_L(\dbf)$ pour all{\'e}ger les notations; nous allons
construire une filtration sur $D_L  = L \otimes_{L_0} D$ {\`a} la
mani{\`e}re de \cite[proposition 5.6]{PC01}.

On a des isomorphismes 
\[  L_n(\!(t)\!) \otimes_L D_L^n = L_n(\!(t)\!) \otimes^{\iota_n}_{\bnrig{,r}{,K}} 
\dbf_r \simeq L_n(\!(t)\!) \otimes^{\iota_n}_{\bnrig{,r}{,K}} 
\mbf_r \] que l'on utilise pour d{\'e}finir 
$\fil^i D_L^n = D_L^n \cap t^i
L_n[\![t]\!] \otimes^{\iota_n}_{\bnrig{,r}{,K}} \mbf_r$ 
et l'isomorphisme $D_L \ra D_L^n$ permet de d{\'e}finir 
$\fil^i D_L$. 
Par d{\'e}finition, l'isomorphisme $D_L \ra D_L^n$ est donn{\'e} par $\mu
\otimes x \mapsto \mu \otimes \iota_n(\phi^n(x))$ et donc
l'isomorphisme $D_L \ra D_L^{n+1}$ co{\"\i}ncide avec la composition 
$D_L \ra D_L^n  \overset{\phi_n}{\ra} D_L^{n+1}$ ce qui fait que
la filtration induite sur $D_L$ ne d{\'e}pend pas de $n$. 

Enfin, on a
$K_n[\![t]\!] \otimes^{\iota_n}_{\bnrig{,r}{,K}} \mbf_r 
= \fil^0 (K_n(\!(t)\!) \otimes_K D_K^n)$  
ce qui fait que $\mbf = \pgm(D) = \pgm(\on{Sol}_L(\dbf))$.
\end{proof}

En rassemblant les th{\'e}or{\`e}mes \ref{eqcat} et \ref{imageessent}, on
trouve:  

\begin{theo}\label{maineqcat}
Le $\otimes$-foncteur exact
$D \mapsto \pgm(D)$, de la cat{\'e}gorie des
$(\phi,N,G_K)$-modules filtr{\'e}s, dans la cat{\'e}gorie des
$(\phi,\Gamma_K)$-modules sur $\bnrig{}{,K}$ dont la connexion
associ{\'e}e est localement triviale, est une {\'e}quivalence de
cat{\'e}gories. 
\end{theo}

En particulier, on a le r{\'e}sultat suivant qui nous servira dans la
suite: 

\begin{coro}\label{sousobj}
Si $D$ est un $(\phi,N,G_{L/K})$-module filtr{\'e}, et
si $\mbf'$ est un sous-$(\phi,\Gamma_K)$-module de $\mbf = \pgm(D)$, alors il
existe un sous-objet $D' \subset D$ tel que $\mbf'=\pgm(D')$.
\end{coro}

\begin{proof}
Par le lemme \ref{soustriv} et la remarque qui le suit, 
la connexion $\nabla_{\mbf}$ est localement triviale sur $\mbf'$ et on peut
donc {\'e}crire $\mbf'=\pgm(D')$ pour un $(\phi,N,G_{M/K})$-module
filtr{\'e} o{\`u} $M/K$ est une extension suffisamment grande. 
L'inclusion $\mbf' \subset \mbf$ nous donne par
fonctorialit{\'e} une inclusion de $(\phi,N,G_{M/K})$-modules
filtr{\'e}s $D' \subset D$, ce qui fait que $G_{M/L}$ agit trivialement
sur $D'$ et donc que $D'$ est un sous-$(\phi,N,G_{L/K})$-module
filtr{\'e} de $D$.
\end{proof}

\section{Pentes de Frobenius}\label{secIV}
Nous calculons dans ce paragraphe les pentes de Frobenius des 
$(\phi,\Gamma_K)$-modules $\pgm(D)$ associ{\'e}s aux
$(\phi,N,G_K)$-modules filtr{\'e}s $D$. En particulier, on montre que $D$
est admissible si et seulement si $\pgm(D)$ est {\'e}tale.

\subsection{Rappels sur les pentes de Frobenius}\label{frobslopes}
Rappelons le th{\'e}or{\`e}me principal (le
th{\'e}or{\`e}me 6.10) de \cite{KK00}. Pour cela, il convient de noter que
les anneaux $\Gamma_{\mathrm{an,con}}$ et $\Gamma_{\mathrm{con}}[1/p]$ de
Kedlaya sont nos anneaux $\bnrig{}{,K}$ et $\bdag{}_K$.
\begin{theo}\label{pentes}
Si $\mbf$ est un $\phi$-module sur $\bnrig{}{,K}$, alors $\mbf$ admet
une filtration canonique 
$\{0\} =\mbf_0 \subset \mbf_1 \subset \cdots \subset \mbf_\ell = \mbf$ 
par des sous $\phi$-modules telle que:
\begin{enumerate}
\item pour $i=1,\cdots,\ell$, le quotient $\mbf_i/\mbf_{i-1}$
est isocline de pente sp{\'e}ciale $s_i$;
\item $s_1 < s_2 < \cdots < s_\ell$;
\item chaque quotient $\mbf_i/\mbf_{i-1}$ peut s'{\'e}crire
$\mbf_i/\mbf_{i-1} = \bnrig{}{,K} \otimes_{\bdag{}_K} \nbf_i$ o{\`u}  
$\nbf_i$ est un $\phi$-module sur $\bdag{}_K$ isocline de pente g{\'e}n{\'e}rique
$s_i$.
\end{enumerate}
De plus, les conditions (1) et (2) ci-dessus d{\'e}terminent la
filtration, et les $\nbf_i$ du (3) sont aussi uniques.
\end{theo}

Nous ne rappelons pas ce que sont les pentes sp{\'e}ciales et
g{\'e}n{\'e}riques, mais signalons que si $\mbf$ est de rang $1$ alors les
pentes sp{\'e}ciales et g{\'e}n{\'e}riques co{\"\i}ncident, et si de plus on
a $\mbf=\bnrig{}{,K} \cdot x$ o{\`u} $\phi(x)=\lambda x$ avec $\lambda \in
\OO_{K_0'}$, alors cette pente est {\'e}gale {\`a} $v_p(\lambda)$. 

Si $\mbf$ est un $(\phi,\Gamma_K)$-module, 
alors, comme l'action de $\Gamma_K$ commute {\`a} $\phi$, 
les $\mbf_i$ sont stables par
$\Gamma_K$ puisque la filtration est 
canonique et les $\nbf_i$ sont
aussi stables par $\Gamma_K$ par unicit{\'e}.

\begin{defi}\label{defetale}
On dit qu'un $\phi$-module $\mbf$ sur $\bnrig{}{,K}$ est {\'e}tale si
dans le th{\'e}or{\`e}me \ref{pentes} ci-dessus, $\ell=1$ et $s_1=0$,
c'est-{\`a}-dire s'il existe un $(\phi,\Gamma_K)$-module {\'e}tale
(i.e. de pente g{\'e}n{\'e}rique nulle)
$\mbf^{\dagger}$ sur $\bdag{}_K$ tel que $\mbf = \bnrig{}{,K}
\otimes_{\bdag{}_K} \mbf^{\dagger}$.
\end{defi}

\begin{prop}\label{souspente}
Si $\mbf$ est un $\phi$-module, dont les pentes $s_1 < s_2 < \cdots <
s_\ell$ sont $\geq 0$, et $x \in \mbf$ est tel que $\phi(x)=\lambda x$
avec $\lambda \in K_0'$, alors $\lambda \in \OO_{K_0'}$. 
\end{prop}

\begin{proof}
Le vecteur $x$ est un vecteur propre de Frobenius. 
La d{\'e}monstration de \cite[theorem 6.10]{KK00}
montre que $\mbf_1 \subset \mbf$ est de pente {\'e}gale {\`a} la plus
petite pente sp{\'e}ciale de $\mbf$. Les pentes sp{\'e}ciales de $\mbf$
sont d{\'e}finies {\`a} partir des vecteurs propres de Frobenius (cf
\cite[\S 4.4]{KK00}) ce qui fait que si $s_1 \geq 0$, alors
$v_p(\lambda) \geq 0$. 
\end{proof}

\subsection{Calcul des pentes de $\pgm(D)$}\label{calcslopes}
Le r{\'e}sultat principal de ce chapitre est le th{\'e}or{\`e}me ci-dessous:
\begin{theo}\label{detslope}
Si $D$ est un $(\phi,N,G_{L/K})$-module filtr{\'e}, alors la pente
de $\det \pgm(D)$ est {\'e}gale {\`a} $t_N(D)-t_H(D)$.  
\end{theo}

\begin{proof}
Comme le foncteur $D \mapsto \pgm(D)$ est un $\otimes$-foncteur exact,
on a $\det \pgm(D) = \pgm(\det D)$ et d'autre part on a par d{\'e}finition
$t_N(D)=t_N(\det D)$ et $t_H(D)=t_H(\det D)$ ce qui fait qu'il suffit
de montrer le th{\'e}or{\`e}me quand $D$ est de rang $1$. Si $e$ est une
base de $D$ telle que $\phi(e)=p^\nu \lambda_0 e$ o{\`u} $\lambda_0 \in
\OO_{K_0}^*$ et $t_H(e)=\eta$, alors on voit que $\pgm(D)=
\bnrig{}{,K}  t^{-\eta} \otimes e$ et que $\phi(t^{-\eta} \otimes e) =
p^{\nu-\eta}\lambda_0 \cdot t^{-\eta} \otimes e$ ce qui fait que la
pente de $\pgm(D)$ est {\'e}gale {\`a} $\nu-\eta$ et vaut bien
$t_N(D)-t_H(D)$. 
\end{proof}

\begin{prop}\label{admiss}
Si $D$ est un $(\phi,N,G_{L/K})$-module filtr{\'e}, alors
$D$ est admissible si et seulement si
$\pgm(D)$ est un $(\phi,\Gamma_K)$-module {\'e}tale.
\end{prop}

\begin{proof}
Supposons tout d'abord que $D$ est admissible.
Le th{\'e}or{\`e}me \cite[theorem 6.10]{KK00} de Kedlaya 
rappel{\'e} ci-dessus montre que $\pgm(D)$
admet une filtration canonique par des sous $(\phi,\Gamma_K)$-modules
isoclines de pentes croissantes. La somme de ces pentes (compt{\'e}es
avec multiplicit{\'e}s) est la pente de $\det \pgm(D)$ (cf \cite[prop
5.13]{KK00}) et vaut donc $t_N(D)-t_H(D) = 0$. Pour montrer que
$\pgm(D)$ est isocline de pente nulle, il suffit donc de montrer que
les pentes de $\pgm(D)$ sont $\geq 0$. Par le corollaire \ref{sousobj},
tout sous-objet de $\pgm(D)$
est de la forme $\pgm(D')$ o{\`u} $D' \subset D$ et la pente de $\det
\pgm(D')$ vaut $t_N(D')-t_H(D') \geq 0$ puisque $D$ est suppos{\'e}
admissible. On en conclut que $\pgm(D)$ ne peut pas contenir de
sous-objet isocline de pente $<0$ et donc que $\pgm(D)$ est
{\'e}tale. 

Montrons maintenant que si $\pgm(D)$ est {\'e}tale, alors $D$ est
admissible. La pente de $\det \pgm(D)$ est nulle et donc
$t_N(D)-t_H(D) = 0$. Si $D'$ est un sous-objet de $D$, de dimension
$d'$, alors $\det(D')$ est de dimension $1$ dans $\wedge^{d'} D$ 
et $\pgm(\det D')$ est un sous-$\phi$-module de rang $1$ de
$\pgm(\wedge^{d'} D)$.
Par la proposition \ref{souspente}, la pente 
de $\pgm(\det D')$ est $\geq 0$ ce qui fait que $t_N(D')-t_H(D') \geq
0$. On en conclut que $D'$ est admissible.
\end{proof}

\begin{rema}\label{thisisfadm}
Comme on le verra au chapitre \ref{secV}, la proposition \ref{admiss}
ci-dessus permet de red{\'e}montrer le th{\'e}or{\`e}me d'admissibilit{\'e} de
Colmez-Fontaine en utilisant la correspondance entre repr{\'e}sentations
$p$-adiques et $(\phi,\Gamma_K)$-modules {\'e}tales.
\end{rema}

Le th{\'e}or{\`e}me suivant est une cons{\'e}quence imm{\'e}diate de la
proposition \ref{admiss} ci-dessus.

\begin{theo}\label{maineqcatadm}
Le $\otimes$-foncteur exact
$D \mapsto \pgm(D)$, de la cat{\'e}gorie des
$(\phi,N,G_K)$-modules filtr{\'e}s admissibles, 
dans la cat{\'e}gorie des $(\phi,\Gamma_K)$-modules 
{\'e}tales sur $\bnrig{}{,K}$ dont la connexion
associ{\'e}e est localement triviale, est une {\'e}quivalence de
cat{\'e}gories. 
\end{theo}

\begin{rema}\label{faltot}
Ce th{\'e}or{\`e}me permet notamment de retrouver le th{\'e}or{\`e}me de
Faltings-Totaro qui affirme que la cat{\'e}gorie des
$(\phi,N,G_K)$-modules filtr{\'e}s admissibles est stable par produits
tensoriels.
\end{rema}

\begin{rema}\label{allslopes}
En g{\'e}n{\'e}ral, on peut se demander comment calculer les pentes de
$\pgm(D)$. Posons $\mu(D)=(t_N(D)-t_H(D))/\dim D$. Si $D$ est
irr{\'e}ductible, alors $\pgm(D)$ est aussi irr{\'e}ductible et 
le th{\'e}or{\`e}me \ref{detslope} montre que $\pgm(D)$ est isocline de
pente $\mu(D)$. Cela sugg{\`e}re un proc{\'e}d{\'e} pour calculer les pentes
de $\pgm(D)$ en g{\'e}n{\'e}ral: parmi tous les sous-objets 
irr{\'e}ductibles $D' \subset
D$, en choisir un $D_{\min}$ qui minimise $\mu(D')$. 
Les pentes de $\pgm(D)$ sont
alors: $\mu(D_{\min})$ et les pentes de $D/D_{\min}$. 
\end{rema}

\begin{rema}\label{nocontradict}
Si $D$ est un $(\phi,N)$-module filtr{\'e} admissible, alors $\pgm(D)$
peut tr{\`e}s bien contenir des sous-objets de pente $>0$. Ceci n'est
pas en contradiction avec le th{\'e}or{\`e}me de Kedlaya, car c'est la
filtration par des pentes croissantes qui est canonique. 
\end{rema}

\begin{exem}
Voici quelques exemples de calcul de $\pgm(D)$ pour des
$\phi$-modules filtr{\'e}s $D$ relatifs {\`a} $\Qp$. 
Dans tous les cas, $D=\Qp e \oplus \Qp f$. 
\begin{enumerate}
\item $\phi(e)=e$, $\phi(f)=pf$, $\fil^0 D=D$, $\fil^1 D = \Qp(e+f)$
  et $\fil^2 D = \{0\}$. Ce $D$ est admissible.
  Une base de $\pgm(D)$ est donn{\'e}e par $e$ et
  $\alpha e +f$ o{\`u} $\alpha$ est une fonction telle que
  $\alpha(\zeta_{p^n}-1)=p^{-n}$ pour $n \gg 0$. Il y a deux
  sous-$\phi$-modules dans $\pgm(D)$: celui engendr{\'e} par $e$ 
  est de pente $0$ et celui engendr{\'e} par $f$ est de pente $1$.
\item  $\phi(e)=e$, $\phi(f)=pf$, $\fil^0 D=D$, $\fil^1 D = \Qp e$
  et $\fil^2 D = \{0\}$. Ce $D$ n'est pas admissible. Une base de
  $\pgm(D)$ est donn{\'e}e par $t^{-1}e$ et $f$. Il y a deux
  sous-$\phi$-modules dans $\pgm(D)$: celui engendr{\'e} par $t^{-1}e$
  est de pente $-1$ et celui engendr{\'e} par $f$ est de pente $1$.
\item $\phi(e)=p^2 f$, $\phi(f)=e$, $\fil^0 D=D$, $\fil^{1,2} D =
  \Qp e$
  et $\fil^3 D = \{0\}$. Ce $D$ est admissible.
  Une base de $\pgm(D)$ est donn{\'e}e par $e/t_+^2$ et $f/t_-^2$ o{\`u}
  les fonctions \[ t_+(X)=\prod_{n=1}^\infty  
  \frac{(1+X)^{p^{2n}}-1}{p((1+X)^{p^{2n-1}}-1)}
  \quad\text{et}\quad t_-(X)=\prod_{n=1}^\infty  
  \frac{(1+X)^{p^{2n-1}}-1}{p((1+X)^{p^{2n-2}}-1)} \]
  sont les produits partiels pairs et
  impairs de $t=\log(1+X)$. Il y
  a deux sous-$\phi$-modules dans $\pgm(D)$, engendr{\'e}s par $e \pm
  pf$, qui sont tous les deux de pente $1$. 
\end{enumerate}
\end{exem}

\section{Applications aux repr{\'e}sentations $p$-adiques}\label{secV}
L'objet de ce chapitre est de donner des applications 
des constructions ci-dessus aux
repr{\'e}sentations $p$-adiques.  Dans tout ce chapitre, une
repr{\'e}sentation $p$-adique de $G_K$ est un $\Qp$-espace 
vectoriel de dimension finie muni d'une 
action lin{\'e}aire et continue
de $G_K$. 

\subsection{Anneaux de Fontaine et repr{\'e}sentations
  $p$-adiques}\label{fontrap} 
Dans ce paragraphe, nous donnons de brefs rappels sur certains aspects
de la th{\'e}orie des repr{\'e}sentations $p$-adiques. Afin d'{\'e}tudier
les repr{\'e}sentations $p$-adiques, on construit un certains
nombres d'anneaux de p{\'e}riodes, et nous avons besoin des anneaux
$\bst$ et $\bdR$ d{\'e}finis dans \cite{Bu88per} et de l'anneau $\bdag{}$
d{\'e}fini dans \cite{CC98}. 
L'anneau $\bst$ est une $\Qp$-alg{\`e}bre qui est aussi un
$(\phi,N,G_K)$-module (mais $\bst$ est de dimension infinie 
et $G_K$ n'agit pas {\`a} travers un quotient 
fini cette fois) et $\bdR$ est un corps qui est aussi une
$\Qp$-alg{\`e}bre filtr{\'e}e. On de plus une injection $K \otimes_{K_0}
\bst \ra \bdR$ ce qui fait que l'on peut voir $\bst$ comme une
sorte de $(\phi,N,G_K)$-module filtr{\'e}. 

{\'E}tant donn{\'e}e une repr{\'e}sentation $p$-adique $V$ de $G_K$, 
le $K_0$-espace vectoriel $(\bst \otimes_{\Qp} V)^{G_K}$ est de
dimension $\leq \dim(V)$ et on dit
que $V$ est semi-stable si sa dimension est $=\dim(V)$. Si $V$ n'est
pas semi-stable mais le devient quand on la restreint {\`a} un
sous-groupe ouvert $G_L$ de $G_K$, alors $\dstp{L}(V) = (\bst
\otimes_{\Qp} V)^{G_L}$ est un $(\phi,N,G_{L/K})$-module filtr{\'e}
admissible. En fait, $L \otimes_{L_0} \dstp{L}(V)$ s'identifie
naturellement {\`a} $L \otimes_K \ddR(V)$ o{\`u} $\ddR(V)=(\bdR
\otimes_{\Qp} V)^{G_K}$. 

Plus g{\'e}n{\'e}ralement, on dit qu'une repr{\'e}sentation $p$-adique $V$
est de de Rham si le $K$-espace vectoriel $(\bdR \otimes_{\Qp}
V)^{G_K}$ est de dimension $\dim(V)$ et la construction pr{\'e}c{\'e}dente
montre que les repr{\'e}sentations potentiellement 
semi-stables sont de de Rham. 
La r{\'e}ciproque est aussi vraie, ainsi que Fontaine l'avait
conjectur{\'e} dans \cite[\S 6.2]{Bu88sst}
(c'est la conjecture de monodromie pour
les repr{\'e}sentations $p$-adiques. 
Voir le {\og s{\'e}minaire Bourbaki \fg} \cite{PC01} pour des d{\'e}tails
sur cette conjecture; la d{\'e}monstration
du fait que la conjecture de Crew (la conjecture de monodromie pour
les {\'e}quations diff{\'e}rentielles $p$-adiques) implique
la conjecture de monodromie pour
les repr{\'e}sentations  $p$-adiques se trouve dans \cite{LB02} et la
conjecture de Crew est d{\'e}montr{\'e}e dans \cite{YA01, KK00, ZM01}. 
Des d{\'e}monstrations {\og directes \fg} de la conjecture de monodromie pour
les repr{\'e}sentations $p$-adiques se trouvent dans \cite{PC03,F04}). 

Dans une autre direction, 
le corps $\bdag{}$ est muni d'un Frobenius $\phi$ et d'une action de
$G_K$ telle que $(\bdag{})^{H_K}=\bdag{}_K$.
Si $\ddag{}$ est un $(\phi,\Gamma_K)$-module
{\'e}tale et de dimension $d$ sur $\bdag{}_K$, 
alors on peut montrer que
$V(\ddag{})=(\bdag{} \otimes_{\bdag{}_K} 
\ddag{})^{\phi=1}$ est un $\Qp$-espace
vectoriel de dimension $d$ qui h{\'e}rite d'une action de $G_K$. Si on
combine les r{\'e}sultats de Fontaine (cf \cite{F90}) et le th{\'e}or{\`e}me
de Cherbonnier-Colmez (cf \cite{CC98}), 
on trouve que le foncteur $\ddag{} \ra
V(\ddag{})$ est une {\'e}quivalence de cat{\'e}gories de la cat{\'e}gorie
des $(\phi,\Gamma_K)$-modules {\'e}tales sur $\bdag{}_K$ vers la
cat{\'e}gorie des repr{\'e}sentations $p$-adiques de $G_K$; on note $V
\mapsto \ddag{}(V)$ l'inverse de ce foncteur.

Si $V$ est de de Rham, alors on peut montrer (cf \cite[proposition 
  3.25]{F00}) que la connexion $\nabla_{\dbf}$ sur $\dnrig{}(V) =
\bnrig{}{,K} \otimes_{\bdag{}_K} \ddag{}(V)$ est localement triviale,
et le th{\'e}or{\`e}me \ref{imageessent} 
(qui dans ce cas est aussi donn{\'e} dans \cite[th{\'e}or{\`e}me
  5.10]{LB02}) montre
qu'il existe alors un unique $(\phi,\Gamma_K)$-module $\ndr(V)$ de rang
$\dim(V)$ contenu dans $\dnrig{}(V)$ et tel que $\nabla_{\dbf}(\ndr(V))
\subset t \ndr(V)$. 
Cette construction est d'ailleurs le point de d{\'e}part d'une 
d{\'e}monstration du th{\'e}or{\`e}me de monodromie pour les
repr{\'e}sentations $p$-adiques.

\subsection{Repr{\'e}sentations potentiellement
  semi-stables}\label{potsstreps} 
Si $V$ est une repr{\'e}sentation $p$-adique de $G_K$ qui devient
semi-stable quand on la restreint {\`a} $G_L$ pour une extension
galoisienne finie $L/K$, alors $\dstp{L}(V) = (\bst
\otimes_{\Qp} V)^{G_L}$ est un $(\phi,N,G_{L/K})$-module filtr{\'e}
admissible. L'objet de ce paragraphe est de montrer comment on peut
utiliser le th{\'e}or{\`e}me \ref{maineqcatadm} pour donner une nouvelle
d{\'e}monstration du th{\'e}or{\`e}me de Colmez-Fontaine (cf
\cite[th{\'e}or{\`e}me A]{CF00}) rappel{\'e} ci-dessous:
\begin{theo}\label{colfont}
Si $D$ est un $(\phi,N,G_{L/K})$-module filtr{\'e}
admissible, alors il existe une repr{\'e}sentation 
$p$-adique $V$ de $G_K$ qui devient
semi-stable quand on la restreint 
{\`a} $G_L$ et telle que $D = \dstp{L}(V)$. 
\end{theo}

\begin{proof}
Si $D$ est un $(\phi,N,G_{L/K})$-module filtr{\'e}
admissible, alors le th{\'e}or{\`e}me \ref{maineqcatadm} montre que
$\pgm(D)$ est un $(\phi,\Gamma_K)$-module {\'e}tale sur $\bnrig{}{,K}$
ce qui fait que l'on peut {\'e}crire $\pgm(D) = \bnrig{}{,K}
\otimes_{\bdag{}_K} \ddag{}$ o{\`u} $\ddag{}$ 
est un $(\phi,\Gamma_K)$-module
{\'e}tale sur $\bdag{}_K$. 
Par la construction $\ddag{} \mapsto V(\ddag{})$
rappel{\'e}e au paragraphe pr{\'e}c{\'e}dent, 
il existe une repr{\'e}sentation 
$p$-adique $V$ de $G_K$ telle que $\ddag{} = \ddag{}(V)$. 

Rappelons que par \cite[th{\'e}or{\`e}me 3.6]{LB02}, si $V$ est une
repr{\'e}sentation $p$-adique de $G_K$, alors $\dstp{L}(V) =
(\bnrig{}{,L}[\ell_X,1/t] 
\otimes_{\bnrig{}{,K}} \dnrig{}(V))^{\Gamma_L}$. Dans notre cas, le
fait que \[ \bnrig{}{,L}[\ell_X,1/t] \otimes_{\bnrig{}{,K}} 
\pgm(D) = \bnrig{}{,L}[\ell_X,1/t] \otimes_{L_0} D \] et que
$(\bnrig{}{,L}[\ell_X,1/t])^{\Gamma_L}=L_0$ montrent que
$D = (\bnrig{}{,L}[\ell_X,1/t] 
\otimes_{\bnrig{}{,K}} \dnrig{}(V))^{\Gamma_L}$ et donc que
$D=\dstp{L}(V)$ en tant que $(\phi,N,G_{L/K})$-modules. Il reste {\`a}
voir que la filtration de $L \otimes_{L_0} D$ 
construite dans le th{\'e}or{\`e}me \ref{imageessent} co{\"\i}ncide avec
la filtration provenant de l'isomorphisme $(L \otimes_{L_0}
\dstp{L}(V))^{G_{L/K}} = \ddR(V)$. Pour cela, il suffit de constater que
par \cite[\S 2.4]{LB02}, 
l'application $\iota_n$ envoie $\bnrig{}{,L}[\ell_X,1/t]
\otimes_{\bnrig{}{,K}} \dnrig{}(V)$
dans $\bdR \otimes_{\Qp} V$ et donc $L_n(\!(t)\!) \otimes_L D_L$ dans
$\bdR \otimes_{\Qp} V$ ce qui fait que, comme pour $n$ assez grand on
a (cf \cite[\S 5.3]{LB02} et \cite[\S 3]{F00}):
\[ L_n[\![t]\!] \otimes^{\iota_n}_{\bnrig{,r}{,K}} \dnrig{,r}(V) =
\on{Fil}^0 (L_n(\!(t)\!) \otimes_K \ddR(V)), \] 
que la filtration construite dans
le th{\'e}or{\`e}me \ref{imageessent} co{\"\i}ncide avec
la filtration provenant de celle de $\bdR \otimes_{\Qp} V$, et donc de
$\ddR(V)$. 

On a donc construit une repr{\'e}sentation $p$-adique $V$ 
de $G_K$ dont la
restriction {\`a} $G_L$ est semi-stable et telle que 
$\dstp{L}(V) = D$ en
tant que $(\phi,N,G_{L/K})$-modules filtr{\'e}s.
\end{proof}

\begin{rema}
La d{\'e}monstration ci-dessus utilise la construction $D \mapsto
\pgm(D)$ mais pas la caract{\'e}risation de l'image essentielle de ce
foncteur, ce qui fait que notre d{\'e}monstration n'utilise pas le
th{\'e}or{\`e}me de monodromie $p$-adique (mais on utilise la filtration
par les pentes de Frobenius).
\end{rema}

Pour finir, nous allons r{\'e}capituler les diff{\'e}rents objets que l'on
associe {\`a} une repr{\'e}sentation $p$-adique $V$ de $G_K$ semi-stable quand
on la restreint {\`a} $G_L$. Ces objets sont: 
\begin{enumerate}
\item le $(\phi,N,G_{L/K})$-module filtr{\'e} $\dstp{L}(V)$;
\item l'{\'e}quation diff{\'e}rentielle $p$-adique $\ndr(V)$;
\item le $(\phi,\Gamma_K)$-module
{\'e}tale sur l'anneau de Robba $\dnrig{}(V)$.
\end{enumerate}

Ces objets sont reli{\'e}s entre eux de la mani{\`e}re suivante (rappelons
que par le th{\'e}or{\`e}me \ref{imageessent}, la donn{\'e}e de
$\dnrig{}(V)$ d{\'e}termine une filtration sur $L \otimes_{L_0}
\on{Sol}_L(\ndr(V))$):  
\begin{theo}\label{recapst}
Si $V$ est une repr{\'e}sentation $p$-adique de $G_K$ semi-stable quand
on la restreint {\`a} $G_L$, alors:
\begin{enumerate}
\item $\dstp{L}(V) = \on{Sol}_L(\ndr(V))$ avec la filtration
provenant de $\dnrig{}(V)$;  
\item $\dnrig{}(V)=\pgm(\dstp{L}(V))$;
\item $\ndr(V)=(\bnrig{}{,L}[\ell_X] \otimes_{L_0} \dstp{L}(V))^{H_K,N=0}$.
\end{enumerate}
Ces identifications sont compatibles {\`a} toutes les structures en pr{\'e}sence.
\end{theo}

\begin{proof}
La d{\'e}monstration du th{\'e}or{\`e}me \ref{colfont} montre que l'on a 
$\dnrig{}(V)=\pgm(\dstp{L}(V))$.  
Si on pose $\nbf=(\bnrig{}{,L}[\ell_X] \otimes_{L_0}
\dstp{L}(V))^{H_K,N=0}$, alors 
$\nbf[1/t] = \dnrig{}(V)[1/t]$ et $\nbf$ est un $(\phi,\Gamma_K)$-module
tel que $\nabla_{\nbf}(\nbf) \subset t \nbf$ ce qui fait que
$\nbf=\ndr(V)$ par \cite[th{\'e}or{\`e}me 5.10]{LB02}. Enfin, le fait que
$\dnrig{}(V)=\pgm(\dstp{L}(V))$ et que
$\ndr(V)[1/t]=\pgm(\dstp{L}(V))[1/t]$ montrent que 
$\dstp{L}(V) = \on{Sol}_L(\ndr(V))$ avec la filtration
provenant de $\dnrig{}(V)$.
\end{proof}

Remarquons que le fait qu'on retrouve l'action de $G_{L/K}$ {\`a}
partir de $\ndr(V)$ avait {\'e}t{\'e} observ{\'e} par Marmora (cf \cite[\S
4.2]{AM03}).  

\subsection{Construction de $\ddag{}(V)$}\label{secord}
Dans tout ce paragraphe, 
$V$ est une repr{\'e}sentation semi-stable de $G_K$. Comme on l'a
signal{\'e} au paragraphe pr{\'e}c{\'e}dent, on peut r{\'e}cup{\'e}rer
$\dnrig{}(V)$ {\`a} partir de $\dst(V)$ par la recette
$\dnrig{}(V)=\pgm(\dst(V))$. Plus explicitement, si on regarde comment
le foncteur $\pgm$ est d{\'e}fini, on voit que $\dnrig{}(V)$ est
l'ensemble des $x \in \bnrig{,r}{,K}[\ell_X,1/t] \otimes_{K_0} \dst(V)$
tels que:
\begin{enumerate}
\item $N(x)=0$;
\item $\phi^{-n}(x) \in \fil^0 (K_n(\!(t)\!) \otimes_K \ddR(V))$ pour tout
  $n \geq n(r)$.
\end{enumerate}
L'objet de ce chapitre est de montrer
comment r{\'e}cup{\'e}rer le $\bdag{,r}_K$-module 
$\ddag{,r}(V)$ {\`a} partir de
$\bnrig{,r}{,K}[\ell_X,1/t] \otimes_{K_0} \dst(V)$.

Rappelons que l'on a construit
dans \cite[\S 2]{LB02} un anneau $\btrig{,r}{}$ qui contient $\bnrig{,r}{,K}$
et aussi l'anneau $\btrigplus{}{}$ de \cite[\S 1.2]{LB02} (dont on
rappelle bri{\`e}vement la construction ci-dessous). 
Nous allons rappeler la d{\'e}finition de l'ordre
d'un {\'e}l{\'e}ment $x \in \btrig{,r}{}$. Pour cela, rappelons 
(cf \cite[\S 2.1]{LB02}) que l'anneau
$\btdag{,r}$ construit dans \cite{CC98} 
est un sous-anneau du corps des fractions $W(\et)[1/p]$
de l'anneau des vecteurs de Witt
sur un corps valu{\'e} $\et$ et que si $x=\sum_{k \gg -\infty}
p^k [x_k] \in \btdag{,r}$ et si $I$ est un intervalle compris dans
$[r;+\infty[$, alors la formule:
\[ V_I(x) = \left\lfloor \inf_{\alpha \in I} \inf_{k \in \ZZ} k +
\frac{p-1}{p\alpha} v_{\et}(x_k) \right\rfloor \]
d{\'e}finit une valuation sur $\btdag{,r}$ et que 
par d{\'e}finition (cf \cite[\S 2.3]{LB02}), $\btrig{,r}{}$ est le
compl{\'e}t{\'e} de $\btdag{,r}$ pour la topologie de Fr{\'e}chet d{\'e}finie
par l'ensemble des valuations $\{ V_{[r;s]} \}_{s \geq r}$. 

Si $s \geq r$, on a par cons{\'e}quent une valuation
$V_{[s;s]}$ sur $\btrig{,r}{}$ dont 
on peut montrer (cf \cite[\S 7.2]{PC03} et \cite[lemme 2.7]{LB02}) que 
la restriction {\`a} $\bnrig{,r}{,K}$
co{\"\i}ncide avec la partie enti{\`e}re de la 
valuation associ{\'e}e {\`a} la norme {\og sup sur
la couronne de rayon $p^{-1/e_ks}$ \fg}. Posons
$\rho=(p-1)/p$. Si $x \in \btrig{,r}{}$ et $n$ est assez grand,
alors $\phi^{-n}(x) \in \btrig{,\rho}{}$. On pose $I(x) = \{ s \in
\RR$ tels que la suite $\{ ns+V_{[\rho,\rho]}(\phi^{-n}(x))\}$ est
born{\'e}e inf{\'e}rieurement$\}$, ce qui fait que soit $I(x)$ est vide,
soit il existe $\alpha \in \RR \cup \{ - \infty \}$ tel que 
$I(x)=[\alpha;\infty[$ ou $I(x)=]\alpha;\infty[$.

\begin{defi}\label{ord}
Si $I(x)=[\alpha;\infty[$, on dit que $x$ est d'ordre
$\alpha$ et si $I(x)=]\alpha;\infty[$ on dit que $x$ est d'ordre
$\alpha^+$. On {\'e}crit alors $\on{ord}(x)=\alpha$ ou
$\on{ord}(x)=\alpha^+$. 
\end{defi}

Rappelons que l'on a {\'e}crit $\ell_X$ pour l'{\'e}l{\'e}ment $\log(\pi)$ 
de \cite[\S 2.4]{LB02} et que $\log(\pi) = \logpi +
\log([\overline{\pi}]/\pi)$ avec $\log([\overline{\pi}]/\pi) \in
\btdag{}_{\Qp}$ ce qui nous permet
de prolonger $\on{ord}(\cdot)$
{\`a} $\btrig{}{}[\ell_X] = \btrig{}{}[\logpi]$ en 
d{\'e}cidant que $\on{ord}(x_0+x_1 \logpi + \cdots+x_k
\logpi^k)=\sup_{0 \leq i \leq k} \on{ord}(x_i)+i$. 

\begin{prop}\label{propord}
La fonction $\on{ord}$ v{\'e}rifie les propri{\'e}t{\'e}s suivantes:
\begin{enumerate}
\item $\on{ord}(x+y)
\leq \max(\on{ord}(x),\on{ord}(y))$ et $\on{ord}(xy)
\leq \on{ord}(x)+\on{ord}(y)$;\label{ordval}
\item $t=\log(1+X)$ est d'ordre $1$ et $\on{ord}(tx)=\on{ord}(x)+1$;\label{ordt}
\item $x \in \btrig{}{}$ appartient {\`a} $\btdag{}$ si et
seulement si $\on{ord}(x) \leq 0$;\label{weier}
\item Si $f(X_K) = \sum_{i \in \ZZ} a_i X_K^i 
\in \bnrig{,r}{,K}$, et $s \geq 0$, 
alors notre d{\'e}finition de l'ordre
co{\"\i}n\-cide avec la d{\'e}finition habituelle, c'est-{\`a}-dire que $s
\in I(x)$ si et seulement si 
$\{ v_p(a_i) +  s \log(i)/\log(p) \}_{i \geq 1}$ est
born{\'e}e inf{\'e}rieurement.\label{ordusual}  
\end{enumerate}
\end{prop}

\begin{proof}
Le point (\ref{ordval}) est imm{\'e}diat. Le point (\ref{ordt}) suit
du fait que $\phi(t)=pt$. Le point (\ref{weier}) suit du fait que  
$x \in \btrig{,r}{}$ appartient {\`a} $\btdag{,r}$ si et
seulement si $\{V_{[r;s]}(x)\}_{s \geq r}$ est born{\'e}e
inf{\'e}rieurement, c'est-{\`a}-dire si $0 \in I(x)$.

Pour montrer le point (\ref{ordusual}), rappelons (cf \cite[\S
2.1]{LB02}) que $V_{[r;s]}(\phi^{-1}(x)) = V_{[pr;ps]}(x)$ et donc
que $s \in I(x)$ si et seulement si la suite 
$\{ ns + V_{[\rho_n,\rho_n]}(x)\}$ est born{\'e}e inf{\'e}rieurement, avec 
$\rho_n=p^{n-1}(p-1)$ (ce qui fait que $\rho=\rho_0$). 
Si $f(X_K) \in \bnrig{,r}{,K}$, et $s \geq r$, 
alors $V_{[s;s]}(f(X_K)) = \lfloor v_p(\sup_{z \in C_K[s;s]}|f(z)|_p) \rfloor$,
et par la th{\'e}orie classique des fonctions holomorphes, on a:
\[ v_p \left(\sup_{z \in C_K[\rho_n;\rho_n]} 
|f(z)|_p \right) = \inf_{i \in \ZZ} \left( v_p(a_i) +
\frac{i}{e_K p^{n-1}(p-1)} \right), \]
ce qui fait que pour montrer le point (\ref{ordusual}), il suffit de
montrer que les deux propri{\'e}t{\'e}s suivantes sont {\'e}quivalentes:
\begin{enumerate}
\item la suite $\{ v_p(a_i) +  s \log(i)/\log(p) \}_{i \geq 1}$ est
born{\'e}e inf{\'e}rieurement;
\item la suite double $\{ ns + v_p(a_i) +
i/(e_K p^{n-1}(p-1)) \}_{i \geq 1, n \geq n(r)}$ est born{\'e}e inf{\'e}rieurement.
\end{enumerate}
Ceci est un exercice d'analyse (r{\'e}elle!).
\end{proof}

La propri{\'e}t{\'e} (\ref{ordt}) 
de la proposition ci-dessus nous permet d'{\'e}tendre
$\on{ord}$ {\`a} $\btrig{}{}[\ell_X,1/t]$ en posant
$\on{ord}(x)=\on{ord}(t^h x)-h$ pour $h \gg 0$. 

Si $\bcris^+$ d{\'e}note l'anneau construit dans \cite{Bu88per},  
tel que $\bst=\bcris^+[\log[\pibar],1/t]$,
alors l'anneau $\btrigplus{}$ est d{\'e}fini par  
$\btrigplus{} = \cap_{n=0}^{\infty} \phi^n(\bcris^+)$ et a les
propri{\'e}t{\'e}s suivantes:
\begin{enumerate}
\item si $V$ est une repr{\'e}sentation $p$-adique, alors
  $\dst(V)=(\btrigplus{}[\log[\pibar],1/t] \otimes_{\Qp} V)^{G_K}$; 
\item $\btrigplus{} \subset \btrig{,r}{}$ et en fait, c'est le
  compl{\'e}t{\'e} pour la topologie de Fr{\'e}chet du sous-anneau de
  $\btdag{,r}$ form{\'e} des {\'e}l{\'e}ments $\sum_{k \gg -\infty} p^k
  [x_k]$ tels que $v_{\et}(x_k) \geq 0$ pour tout $k$. 
\end{enumerate}

Le lemme suivant est imm{\'e}diat et g{\'e}n{\'e}ralise le (\ref{ordt}) de la
proposition ci-dessus.

\begin{lemm}\label{ordper}
Soit $V$ une repr{\'e}sentation semi-stable de $G_K$ et $D$ un
sous-$\phi$-module de $\dst(V)$ de pente $\alpha \in \QQ$. Si $M = (m_{j,i}) 
\in \on{M}(d,\btrigplus{}[\log[\pibar],1/t])$ est la
matrice d'{\'e}l{\'e}ments de $D$ dans une base de $V$, alors
$\on{ord}(m_{j,i}) = \alpha$. 
\end{lemm}

Nous pouvons maintenant {\'e}noncer et d{\'e}montrer 
le r{\'e}sultat principal de ce chapitre.

\begin{theo}\label{ddag}
Si $V$ est une repr{\'e}sentation semi-stable de $G_K$ et si 
\begin{itemize}
\item $\{e_i\}_{i=1 \cdots d}$ est une base de $\dst(V)$ adapt{\'e}e
{\`a} la d{\'e}composition par les pentes de Frobenius;
\item $N(e_i)=\sum_{j=1}^d n_{j,i} e_j$;
\item $\{f_j\}_{j=1 \cdots d}$ est une base de $\ddR(V)$ adapt{\'e}e {\`a} la filtration
et $\phi^{-n}(e_i) = \sum_{j=1}^d p_{j,i}^{(n)} f_j$, 
\end{itemize}
alors $x = \sum_{i=1}^d x_i \otimes e_i 
\in \bnrig{,r}{,K}[\ell_X,1/t] \otimes_{K_0} \dst(V)$ appartient {\`a}
$\ddag{,r}(V)$ si et seulement si:
\begin{enumerate}
\item $N(x_j)+\sum_{i=1}^d n_{j,i} x_i = 0$ pour $j=1 \cdots d$;
\item $\sum_{i=1}^d \iota_n(x_i)p_{j,i}^{(n)} \in t^{-t_H(f_j)}
  K_n[\![t]\!]$ pour $j=1 \cdots d$ et $n \geq n(r)$;
\item $\on{ord}(x_i) \leq - \on{pente}(e_i)$ pour $i=1 \cdots d$.
\end{enumerate}
\end{theo}

\begin{proof}
Un petit calcul montre que la condition (1) est {\'e}quivalente {\`a}
$N(x)=0$ et que la condition (2) est {\'e}quivalente {\`a} 
$\phi^{-n}(x) \in \fil^0 (K_n(\!(t)\!) \otimes_K \ddR(V))$ pour tout
$n \geq n(r)$, ce qui fait que, comme on l'a rappel{\'e} plus haut, $x
\in \dnrig{,r}(V)$ si et seulement s'il satisfait (1) et (2). 

Supposons donc que  $x \in \dnrig{,r}(V)$
satisfait (3). Si $\{v_i\}_{i=1 \cdots d}$ est une base de $V$ et si
l'on {\'e}crit $e_i = \sum_{j=1}^d m_{j,i}v_j$, alors 
$m_{j,i}  \in \btrigplus{}[\log[\pibar],1/t]$ et le fait que $e_i$ est dans un
$\phi$-module de pente $\on{pente}(e_i)$ implique 
par le lemme \ref{ordper} que
$\on{ord}(m_{j,i}) \leq \on{pente}(e_i)$, ce qui fait que 
si $x=\sum_{i=1}^d x_i e_i = \sum_{j=1}^d y_j
v_j$, alors $\on{ord}(y_j) = \on{ord}(\sum_{i=1}^d x_i m_{j,i}) \leq 0$.  
Par le (\ref{weier}) de la proposition \ref{propord}, cela 
implique que $y_j \in
\btdag{}$ et donc que $y_j \in \bnrig{,r}{} \cap \btdag{} = \bdag{,r}$
ce qui fait que $x \in \ddag{,r}(V)$. La r{\'e}ciproque est imm{\'e}diate. 
\end{proof}

\begin{rema}\label{annul}
Si $K=K_0$ et $N=0$, alors
la condition (2) {\'e}quivaut {\`a}: {\og pour tout $n \geq n(r)$,
la fonction $\sum_{i=1}^d \phi^n(p_{j,i}^{(n)}) x_i(X)$ a un z{\'e}ro
d'ordre au moins $-t_H(f_j)$ en $\eps^{(n)}-1$ \fg}. La stabilit\'e de
$\ddag{,r}(V)$ sous l'action de $\Gamma_K$ montre que l'on peut
remplacer $\eps^{(n)}-1$ par n'importe quelle racine primitive
$p^n$-i\`eme de $1$.
\end{rema}

\appendix

\section{Liste des notations}
Voici une liste des principales notations dans l'ordre o{\`u} elles
apparaissent:

\ref{secI}: $K$, $k_K$, $K_n$, $K_\infty$, $K_0$, $K_0'$, $G_K$,
$H_K$, $\Gamma_K$, $\sigma$.

\ref{pnfildef}: $L$, $G_{L/K}$, $D_L$, $t_N(D)$, $t_H(D)$.

\ref{robbarap}: $F$, $\bdag{,r}_F$, $\phi$, $\bdag{}_F$, $\bb_F$,
$e_K$, $\bb_K$, $\bdag{}_K$, $r_0(K)$, $\bdag{,r}_K$, $C_K[r;s]$,
$\bnrig{,r}{,K}$, $t$, $\bnrig{}{,K}$, $r(K)$, $\iota_n$, $n(r)$, 
$\nabla$, $\partial$, $q$.

\ref{pgmodrap}: $\dbf_r$, $r(\dbf)$, $\otimes^{\iota_n}$.

\ref{recres}: $\phi_n$, $M_n$.

\ref{pntopgm}: $\ell_X$, $\xi_n$, $D_K^n$, $M_n(D)$, $\pgm(D)$.

\ref{liealg}: $\nabla_{\dbf}$, $V_{[r;s]}$.

\ref{eqdiff}: $\partial_{\dbf}$, $S_L(\dbf)$, $\on{Sol}_L(\dbf)$.

\ref{fontrap}: $\bst$, $\bdR$, $\bdag{}$, $\dstp{L}(V)$, $\ddR(V)$,
$\ddag{}(V)$, $V(\ddag{})$, $\dnrig{}(V)$, $\ndr(V)$. 

\ref{secord}: $\btrig{,r}{}$, $\btrigplus{}$, $\btdag{,r}$, $V_{[r;s]}$, 
$\on{ord}$, $\log[\pibar]$, $\bcris^+$.

\section{Erratum {\`a} \cite{LB02}}
Comme cet article fait assez naturellement suite {\`a} 
{\og Repr{\'e}sentations $p$-adiques et {\'e}quations diff{\'e}rentielles
\fg} (\cite{LB02}), il me semble utile de donner un erratum. Je
remercie P. Colmez, J-M. Fontaine, J. Teitelbaum et H. Zhu pour leurs
remarques. 

\begin{itemize}
\item[\textit{Example 2.8, 1:}]
remplacer $\amax^+$ par $\amax$.

\item[\textit{Sections 3.3, 5.5:}]
Kedlaya a compl{\`e}tement modifi{\'e} son article [34] et la plupart des
num{\'e}ros des r{\'e}f{\'e}rences sont donc incorrects.

\item[\textit{Th{\'e}or{\`e}me 4.10:}]
le th{\'e}or{\`e}me 4.10 est en fait d{\^u} {\`a} Forster, voir:
O. Forster, Zur Theorie der Steinschen Algebren und Moduln, Math.
Zeitschrift, 97, p. 376ff, 1967.

\item[\textit{Proposition 2.24:}]
l'application $\log$ n'est bien s{\^u}r pas d{\'e}finie pour $x=0$. De
plus je ne l'ai d{\'e}finie que pour $\atplus$ mais plus tard, je
l'utilise sur $\atdag{}$ (par exemple: $\log(\pi_K)$). Il faut donc
l'{\'e}tendre {\`a} $\atdag{}$ ce que fait Colmez dans \cite{PC03}. On
peut aussi le faire {\og {\`a} la main \fg}.

\item[\textit{D{\'e}monstration du lemme 5.27:}]
remplacer 
$\on{GL}_d(\adag{,r},K)$ par 
$\on{GL}_d(\adag{,r}_K)$ et de m{\^e}me, remplacer 
$\on{M}_d(\adag{,r},K)$ 
par $\on{M}_d(\adag{,r}_K)$.

\item[\textit{Matrices:}]
j'ai la mauvaise habitude d'{\'e}crire les matrices {\og {\`a} l'envers
\fg}, par exemple si $f$ et $g$ sont deux applications
semi-lin{\'e}aires, alors dans mes notations 
$\on{Mat}(fg) = f(\on{Mat}(g))\on{Mat}(f)$. 
Pour retrouver la notation habituelle, il faut tout transposer (ce que
j'ai fait dans mes autres articles).

\item[\textit{D{\'e}monstration de la proposition 5.15:}]
il n'est pas vrai que $\iota_n(N_s)=K_n[\![t]\!] \otimes_K \ddR(V)$. 
Ce qui est vrai, c'est que l'image de $\iota_n$ est dense pour la
topologie $t$-adique. C'est ce qui est d{\'e}montr{\'e} et utilis{\'e} dans
le reste de la preuve.

\item[\textit{p.229 l.3:}]
remplacer $\atplus$ par $\widetilde{\mathbf{A}}$.

\item[\textit{L'anneau $\bdag{}_K$:}]
il est affirm{\'e} que l'anneau $\bdag{}_K$ est un anneau de s{\'e}ries
formelles {\`a} coefficients dans $F$ ce qui n'est pas toujours le
cas. C'est un anneau de s{\'e}ries formelles {\`a} coefficients dans
l'extension maximale non-ramifi{\'e}e de $F$ dans $K_\infty$, qui peut
{\^e}tre plus grande que $F$. 
Comme il est vrai que $(\bdag{}_K)^{\Gamma_K}=F$, cela n'affecte pas
les r{\'e}sultats de l'article, et les d{\'e}monstrations sont presque
inchang{\'e}es. 

\item[\textit{Monodromie:}]
pour retrouver le $(\phi,N)$-module $\dst(\cdot)$, il faut prendre
$N(\log(\pi))=-p/(p-1)$ au lieu de $N(\log(\pi))=-1$.

\item[\textit{Diagramme p. 271:}]
dans le diagramme en haut de la page, remplacer $\nabla_M$ par
la connexion associ{\'e}e {\`a} $\partial_M$.

\item[\textit{Lemme 2.7:}]
remplacer $k \gg 0$ par $k \gg -\infty$ dans $\sum_{k \gg 0} p^k [x_k]$.
\end{itemize}

\end{document}